\documentclass[graybox]{svmult}
\usepackage{amssymb}
\usepackage{amsmath}
\usepackage{amsfonts} 
\usepackage{mathptmx}       
\usepackage{helvet}         
\usepackage{courier}        
\usepackage{type1cm}        
\usepackage{makeidx}         
\usepackage{multicol}        
\usepackage[bottom]{footmisc}

\newcommand\bel[1]{\begin{equation}\label{#1}}

\makeindex     
 
\DeclareMathOperator	\sgn  {sgn}

\newcommand \TT		{\mathcal{T}}
\newcommand \Tcal	{\TT}
\newcommand{\dK}{\partial K}

\newcommand{\sumk} {\sum_{K\in\TT^h}}

\newcommand {\sumkez} {\sum_{\substack{ K \in\Tcal^h\\ \ekz\in\dKz}}}

\newcommand{\Hkez}{H_{K,e^{0}}}

\newcommand{\qq}{{{q}_{K,\ekz}}}

\newcommand{\ukp}{u_{K}^{+}}

\newcommand{\ukm}{u_{K}^-}
\newcommand{\ukezm}{u_{K_{e^0}}^-}
\newcommand{\tu}{\widetilde{u}_{K,e^0}^+}

\newcommand \la 		\langle
\newcommand \ra 		\rangle
\newcommand \be   	{\begin{equation}}
\newcommand \ee   	{\end{equation}}
\newcommand \RR      {\mathbb{R}}
\newcommand \eps     \epsilon

\newcommand \om       {\omega}
\newcommand \Om      {\Omega}
\newcommand \al \alpha 
\newcommand \var      {\varphi}

\newcommand{\ekp}{e_{K}^{+}}
\newcommand{\ekm}{e_{K}^{-}}

\newcommand{\ekz}{e^{0}}

\newcommand{\dKz}{\del^{0}K}
\newcommand{\vep}{ \varphi_{e_K^+}}
\newcommand{\vem}{ \varphi_{e_K^-}}
\newcommand{\vepo}{ \varphi_{e_K^+}^{\omega}}

\newcommand{\vepO}{ \varphi_{e_K^+}^{\Omega}}
\newcommand{\vemO}{ \varphi_{e_K^-}^{\Omega}}

\newcommand{\veO}{ \varphi_{e}^{\Omega}}

\newcommand{\varbO} {\var_{\ekp}^\Omega}
\newcommand{\QQ}{{Q}_{K,\ekz}}
\newcommand \Omegabf {\mbox{\boldmath $\Omega$}}

\newcommand{\sumK}{\sum_{e^{0}\in\del^{0}K}}

\newcommand \iht { i^*_{H_t}}

\newcommand \ihz {i^*_{H_0}}

\def\Xint#1{\mathchoice
{\XXint\displaystyle\textstyle{#1}}%
{\XXint\textstyle\scriptstyle{#1}}%
{\XXint\scriptstyle\scriptscriptstyle{#1}}%
{\XXint\scriptscriptstyle\scriptscriptstyle{#1}}%
\!\int}
\def\XXint#1#2#3{{\setbox0=\hbox{$#1{#2#3}{\int}$}
\vcenter{\hbox{$#2#3$}}\kern-.5\wd0}}

\def\dashint{\Xint\diagup}
\newcommand \bei {\begin{itemize}}
\newcommand \eei {\end{itemize}}
\newcommand \Ucal {\mathcal U}
\newcommand \Qbf {\mathbf Q}
\newcommand \ub {\overline u}
\newcommand \vb {\overline v}
\newcommand \omegab {{\widetilde \omega}}
\newcommand \ubar   	{{\overline u}}
\newcommand \barc {\underline c}
\newcommand \cbar {\overline c}
\newcommand \dive 	{\mbox{div}}

\newcommand \del		\partial
\newcommand \Ecal {\mathcal E}

\newcommand{\auth}{\textsc}
\newcommand \R     {\mathbb{R}}

\newcommand \MM     {M}

\newcommand \vbar   {{\overline v}}

\newcommand \Dcal   {{\mathcal D}}

\newcommand \Rn     	{\RR^d}
\newcommand \Ubar {\widetilde U}
\newcommand \demi   {{1/2}}
\newcommand \alphab {\underline \alpha}
\newcommand \ZZ     {\mathbb{Z}}   
\newcommand \im {\text{Im }}    

\newcommand \RN    {\mathbb{R}^N} 
\newcommand \Ocal   {\mathcal O}    
\newcommand \Mcal {\mathcal M}
\newcommand \Lcal {\mathcal L}

\newcommand \Abar  {{\overline A}}

\newcommand \Rcal {\mathcal R}

\begin{document}

\title*{Structure-preserving shock-capturing methods: late-time asymptotics,  curved geometry, small-scale dissipation, and nonconservative products} 

\titlerunning{Structure-preserving shock-capturing methods} 

\author{Philippe G. LeFloch}

\institute{Laboratoire Jacques-Louis Lions, Centre National de la Recherche Scientifique, Universit\'e Pierre et Marie Curie, 4 Place Jussieu, 75252 Paris, France. \email{contact@philippelefloch.org}
\newline 
To appear in: {\em  
``Advances in Numerical Simulation in Physics and Engineering'',  XIV Spanish-French School Jacques-Louis Lions, 
C. V\'azquez-Cend\'on, F. Coquel, C. Par\'es Editors, Springer Verlag, 2014.}
Final version completed in January 2014.
}


\maketitle



\section{Introduction} 

\subsection*{Objective} 

We present some recent developments on shock capturing methods  
for nonlinear hyperbolic systems of balance laws, whose prototype is the Euler system of compressible fluid flows, 
and especially discuss {structure-preserving} techniques. 
The problems under consideration arise with complex fluids in realistic applications when friction terms, geometrical terms, viscosity and capillarity effects, etc., need to be taken into account in order to achieve a proper description of the physical phenomena.  
For these problems,  it is necessary to design numerical methods that are not only  consistent with the given partial differential equations, but      remain accurate and robust in certain {asymptotic regimes} of physical interest. That is, certain structural properties of these hyperbolic problems (conservation or balance law, equilibrium state, monotonicity properties, etc.) are essential in many applications, and one seeks that the numerical solutions preserve these properties, which is often a very challenging task. 

To be able to design structure-preserving methods, a theoretical analysis of the hyperbolic problems under consideration  must be performed first by investigating certain singular limits as well as certain classes of solutions of physical relevance. The mathematical analysis allows one to exhibit the key properties of solutions and derive {effective equations} that describe the limiting behavior of solutions, etc. This step requires a deep understanding of the initial value problem, as is for instance the case of {small-scale dissipation sensitive,} viscosity-capillarity driven shock waves which, as it turns out, do not satisfy standard entropy criteria; see LeFloch \cite{PLF-book} for a review. Such a study is      in many physical applications involving hyperbolic systems in {nonconservative form,} in order to avoid the appearance of spurious solutions with wrong speed; see Hou and LeFloch \cite{HL0}.

The design of structure-preserving schemes forces us to go beyond the basic property of consistency with the conservative form of the equations, and requires to revisit the standard strategies, based on finite volumes,  finite differences, Runge-Kutta techniques, etc. By mimicking the theoretical analysis at the discrete numerical level, we can arrive at {structure-preserving schemes,} which preserve the relevant structure of the systems and the asymptotic behavior of solutions.

The techniques developed for model problems provide us with the proper tools to takkle the full problems of physical interest. A variety of nonlinear hyperbolic problems arsing in the applications do involve small scales or enjoy important structural or asymptotic properties. 
By going beyond the consistency with the conservation form of the equations, one can now develop a variety of numerical methods that preserve these properties at the discrete level.
By avoid physically wrong solutions, one can understand first the physical phenomena in simplified situations, and next contribute to validate the ``full'' physical models.

We will only review here two techniques which allows one to preserve late-asymptotics and geometrical terms
and, for further reading on this broad topic, we refer to the textbooks \cite{Bouchut, PLF-book, LeVeque-book}, as well as the lecture notes \cite{PLF-Freiburg,PLF-Olso,Acta}.  Another challenging application arises in continuum physics in the regime of {(small) viscosity and capillarity}, which may still drive the propagation of certain (nonclassical undercompressive) shock waves. This is is relevant in material science for the modeling of smart (martensite) materials, as well as in fluid dynamics for the modeling of multiphase flows (for instance in the context of nuclear plants) and for the coupling of physical models across interfaces.  


\subsection*{Preserving late-time asymptotics with stiff relaxation}

In Section~\ref{sec:2}, this strategy is developed for a class of hyperbolic systems with {stiff relaxation in the regime of late times.} Such systems arise in the modeling of a complex multi-fluid flow when two (or more) scales drive the behavior of the flow. 
Many examples from continuum physics fall into the proposed framework, for instance the Euler equations with (possibly nonlinear) 
{friction.}   In performing a singular analysis of these hyperbolic systems, we keep in mind the analogy with the passage from Boltzmann equation (microscopic description) to the Navier-Stokes equations (macroscopic description). Our aim here is, first, to derive via a formal Chapman-Enskog expansion an {effective system of parabolic type} and, second, to design a scheme which provides consistent and accurate discretizations for {all times,} including asymptotically late times. 

Indeed, we propose and analyze a class of {asymptotic-preserving finite volume methods}, which 
are consistent with, both, the given nonlinear hyperbolic system and the effective parabolic system. It thus preserves the late-time asymptotic regime and, importantly, requires only a classical CFL (Courant, Friedrichs, Lewy) condition of hyperbolic  type, rather than a more restrictive, parabolic-type stability condition.  This section is based on the joint work  \cite{BLT,BLR}.


\subsection*{Geometry-preserving finite volume methods}

The second topic of interest here is provided by the class of {hyperbolic conservation laws posed on a curved space.} Such equations are relevant in geophysical applications, for which the prototype is given by shallow water equations on the {sphere with topography.} 
Computations of large-scale atmospheric flows and oceanic motions (involving the Coriolis force, Rosby waves, etc.) requires robust numerical methods. Another motivation is provided conservation laws on moving surfaces describing combusion phenomena. 
We should      astrophysical applications, involving fluids or plasmas, and the study of the propagation of linear waves (wave operator, Dirac equations, etc.) on curved backgrounds of general relativity (such as Schwarschild or, more generally, Kerr spacetime).
These applications provide important examples where the partial differential equations
of interest are naturally posed on a curved manifold.

Scalar conservation laws yield a drastically simplified, yet very challenging, mathematical
model for understanding nonlinear aspects of shock wave propagation on manifolds.
In Section~\ref{sec:3}, based on the work \cite{LO1}, we introduce the {geometry-preserving finite volume method} for hyperbolic balance laws 
formulated on surfaces or, more generally, {manifolds.} 
First, we present some theoretical tools to handle the interplay between the nonlinear waves propagating on solutions 
and the underlying geometry of the problem. A generalization of the standard Kruzkov theory is obtained on a manifold, by 
formulating the hyperbolic equation under consideration from a field of differential forms. 
The proposed finite volume method is geometry-consistent and relies on a coordinate-independent formulation.   
The actual implementation of this finite volume scheme on the sphere is realized in \cite{BLM,BFL1}.


\section{Late-time asymptotics with stiff relaxation}
\label{sec:2}

\subsection{A class of nonlinear hyperbolic systems of balance laws} 
\label{sec:21} 

Consider the following system of partial differential equations  
\bel{100}
\eps \, \del_t U + \del_x F(U) = -{R(U) \over \eps}, \qquad  U=U(t,x) \in \Omega \subset \RN, 
\ee
in which $t >0, \, x \in \RR$ denote the time and space variables and the flux $F: \Omega \to \RN$ is defined on the convex and open subset $\Omega$. The first-order part of \eqref{100} is assumed to be hyperbolic in the sense that the matrix-valued map  
 $A(U) := D_UF(U)$ admits real eigenvalues  and a full basis of eigenvectors.

In order to analyze the singular limit $\eps \to  0$ of late-time and stiff relaxation, we distinguish between two distinct regimes.  
In the hyperbolic-to-hyperbolic regime, one replaces $\eps \del_t U$ by $\del_t U$ and establishes that solutions to 
$$ 
\del_t U + \del_x F(U) = - {R(U) \over \eps}, \qquad  U=U(t,x), 
$$
are driven by an effective system of {hyperbolic type.} 
Such a study was pioneered by Chen, Levermore, and Liu \cite{CLL}. 
On the other hand, in the hyperbolic-to-parabolic regime which is under consideration in the present work, we 
obtain effective equations of {parabolic type.} In the earlier papers \cite{Marcati1,Marcati3}, Marcati et al.~established rigorous 
convergence theorems for several classes of models.  Our objective here is to introduce a general framework to design numercal methods for such problems. 

We make the following assumptions. 

\ 

\noindent{\bf  Condition 1.} There exists an $n \times N$ matrix $Q$ with (maximal) rank $n<N$ such that
\be
Q R(U) = 0, \qquad U \in \Omega,
\ee
hence, $QU \in Q\Omega =: \omega $ satisfies
\be
\eps \, \del_t \big( QU  \big)  + \del_x \big( QF(U) \big) = 0. 
\ee
 
\

\noindent{{\bf  Condition 2.} There exists a map $\Ecal:\omega \subset \RR^N \to\Omega$ describing the equilibria $u \in \omega$, with   
\be
R(\Ecal(u))=0, \qquad  u = Q \, \Ecal(u).
\ee
We      introduce the equilibrium submanifold $\Mcal:= \big\{ U= \Ecal(u) \big\}$. 

\

\noindent{{\bf  Condition 3.}  It is      assumed that 
\be
QF(\Ecal(u))= 0, \qquad  u\in\omega.
\ee
Observe that the term 
$\del_x \big( QF(\Ecal(u)) \big)$ must vanish identitically, so that $QF(\Ecal(u))$ must be a constant, which we normalize to be $0$. 
 
\

\noindent{{\bf  Condition 4.} For all $u\in \omega$, we impose 
\be
\aligned
\dim\Big( \ker(B(\Ecal(u)))\Big)&=n,\\
\ker\big( B(\Ecal(u)) \big) \cap \im\big( B(\Ecal(u))\big) & =\{0\},
\endaligned
\ee
hence, the $N\times N$ matrix $B:=DR_U$ has ``maximal'' kernel on the equilibrium manifold. 
 

\subsection{Models arising in compressible fluid dynamics} 

\subsubsection*{Stiff friction}

We beign with the Euler system for compressible fluids with friction: 
\be
\label{Euler1} 
\begin{aligned}
& \eps \, \del_t\rho +\del_x( \rho v ) =0,
\\
& \eps \, \del_t( \rho v) +\del_x \big( \rho v^2 + p(\rho) \big) = - {\rho v \over \eps}. 
\end{aligned}
\ee
The density $\rho \geq 0$ and the velocity $v$ are the main unknowns, while
the pressure $p:\RR^+\to\RR^+$ is a prescribed function satisfying the hyperbolicity condition $p'(\rho)>0$ (for $\rho>0$). The first-order homogeneous system is strictly hyperbolic and \eqref{Euler1} fits into our
late-time/stiff relaxation framework in Section~\ref{sec:21} if we set  
$$
U=\left(\begin{array}{c}
\rho \\ \rho v
\end{array}\right),
\quad
F(U)=\left(\begin{array}{c}
\rho v \\ \rho v^2 + p(\rho)
\end{array}\right),
\quad
R(U)=\left(\begin{array}{c}
0 \\ \rho v
\end{array}\right)
$$
and
$
Q=(1~0). 
$
The local equilibria $u=\rho$ are found to be scalar-valued with 
$\Ecal(u) = \left(\rho, \, 0\right)^T$
and we immediately check that $QF(\Ecal(u))=0$.


\subsubsection*{Stiff radiative transfer}

The following model arises in the theory of radiative transfer:
\bel{M11} 
\begin{aligned}
& \eps \, \del_t e +\del_x f =\frac{\tau^4-e}{\eps},
\\
& \eps \, \del_t f+\del_x \Big(  \chi\left(f/e\right)e \Big) =
-\frac{f}{\eps},
\\
& \eps \, \del_t\tau=\frac{e-\tau^4}{\eps}. 
\end{aligned}
\ee
The radiative energy $e>0$ and the radiative flux $f$ are the main unknowns, 
restricted so that $| f/e | \le 1$, 
while $\tau>0$ is the temperature. 
The so-called Eddington factor $\chi:[-1,1]\to\RR^+$ is, typically, taken to be   
$\chi(\xi) =\frac{ 3+4\xi^2}{5+2\sqrt{4-3\xi^2}}.$
Again, this system fits within our general framework.


\subsubsection*{Coupling stiff friction and stiff radiative transfert}

By combining the previous two examples together, one can consider to the following coupled Euler/$M1$ model
\bel{EM1} 
\begin{aligned}
& \eps\del_t\rho +\del_x \big( \rho v \big) = 0,
\\
& \eps\del_t\rho v+\del_x \big( \rho v^2+p(\rho) \big) 
= -\frac{\kappa}{\eps} \rho v+\frac{\sigma}{\eps} f,
\\
& \eps\del_t e +\del_x f =0,
\\
& \eps\del_t f+\del_x \Big(  \chi\left(\frac{f}{e}\right)e\Big)
 = -\frac{\sigma}{\eps}f.
\end{aligned}
\ee
Here, $\kappa$ and $\sigma$ are positive constants and, in the applications, a typical choice for the pressure
is 
$p(\rho)=C_p \rho^\eta$
with $C_p\ll 1$ and $\eta>1.$ Now, we should set 
$$
U=\begin{pmatrix}\rho\\ \rho v\\ e\\f \end{pmatrix},
\qquad
F(U)=\begin{pmatrix}\rho v\\ \rho v^2+p(\rho)\\ f\\ \chi(\frac{f}{e})e \end{pmatrix},
\qquad
R(U)=\begin{pmatrix}0\\ \kappa \rho v-\sigma f\\ 0\\ \sigma f \end{pmatrix}, 
$$
and the local equilibria read 
$$
 \Ecal(u)=\begin{pmatrix}\rho\\ 0\\ e\\0 \end{pmatrix},
 \qquad
u=QU=\begin{pmatrix}\rho\\ e\end{pmatrix},
\qquad
Q=\begin{pmatrix}1&0&0&0\\0&0&1&0\end{pmatrix}, 
$$
so that, once again, $QF(\Ecal(u))=0$. 


\subsection{An expansion near equilibria} 

Our singular analysis proceeds with a Chapman-Engskog expansion around a local equilibria $u= u(t,x) \in \omega$. We set 
$$
U^\eps = \Ecal(u) + \eps \, U_1 +\eps^2 \, U_2 + \ldots, 
\qquad 
\quad 
 u:= QU^\eps,  
$$
and requires that 
$
\eps \, \del_t U^\eps + \del_x F(U^\eps) = - R(U^\eps)/\eps. 
$
We thus obtain 
$
QU_1= QU_2 = \ldots = 0 
$
and then 
$$
\aligned
F(U^\eps) & = F(\Ecal(u))  + \eps \, A(\Ecal(u))\,  U_1
+ 
\Ocal(\eps^2),
\\
\frac{R(U^\eps)}{\eps}
& = B(\Ecal(u))\,  U_1 +
   \frac{\eps}{2}D^2_U R(\Ecal(u)).(U_1,U_1) 
 +
   \eps B(\Ecal(u))\,  U_2 + \Ocal(\eps^2).
\endaligned
$$
 In turn, we deduce that 
$$
\begin{aligned}
& \eps \, \del_t \big( \Ecal(u) \big)
+ \del_x \big(  F(\Ecal(u)) \big) + \eps \, \del_x \Big( A(\Ecal(u))\,  U_1 \Big) 
\\ 
&
= -B(\Ecal(u))\,  U_1 -\frac{\eps}{2} D^2_U R(\Ecal(u)).(U_1,U_1)
  -\eps B(\Ecal(u))\,  U_2 +\Ocal(\eps^2).
\end{aligned}
$$
 
The {zero-order terms} imply that $U_1 \in \RN$ satisfies the algebraic system
$$
B(\Ecal(u))\,  U_1 = - \del_x \big( F(\Ecal(u))\big) \in \RN, 
$$ 
which we can solve in $U_1$. At this juncture, we rely on the condition $QU_1 =0$ and the following lemma. 
 
\begin{lemma}[Technical lemma] 
If $C$ is an $N \times N$ matrix satisfying 
$\dim \ker C = n$ and $\ker C \cap \im C = \big\{ 0 \big\}$, and if $Q$ is an
$n \times N$ matrix of rank $n$, then for all $J\in\RR^N$, there exists a unique solution $V\in\RR^N$ to 
$C\,  V =J$ and $QV =0$ 
$QJ=0$.
\end{lemma}

\begin{proposition}[First-order corrector problem]
The first-order term $U_1$ is characterized by  
$B(\Ecal(u))\,  U_1 = -\del_x\big( F(\Ecal(u)) \big)$
and $QU_1 = 0$. 
\end{proposition}

Considering next the {first-order terms}, we arrive at 
$$
\del_t \big( \Ecal(u) \big) + \del_x \Big( A(\Ecal(u))\,  U_1 \Big)
= -\frac{1}{2} D^2_U R(\Ecal(u)).(U_1,U_1)- B(\Ecal(u))\,  U_2
$$
and, after multiplication by $Q$ and using $Q\Ecal(u)=u$, 
$$
\del_t u  + \del_x \Big( Q \, A(\Ecal(u))\,  U_1 \Big)
= -\frac{1}{2} Q \, D^2_U R(\Ecal(u)).(U_1,U_1) - Q \, B(\Ecal(u))\,  U_2.
$$
On the other hand, by differentiating $Q R(U)=0$, we get  
$
 Q \, D^2_U R . (U_1,U_1) \equiv 0$ and $Q \, B \,  U_2 \equiv 0$.  This leads us to the following conclusion. 

\begin{theorem}[Late time/stiff relaxation effective equations] The effective system reads 
$$
\del_t u = -\del_x \Big( Q A(\Ecal(u))\, U_1 \Big) =: \del_x \left(\Mcal(u) \del_x u \right) 
$$
for some $n\times n$ matrix $\Mcal(u)$ and with $U_1$ being the unique solution to 
$$
\begin{aligned}
B (\Ecal(u))\,  U_1 &= - A(\Ecal(u)) \del_x\big( \Ecal(u) \big),
\qquad
QU_1  = 0.
\end{aligned}
$$
\end{theorem}


\subsection{Mathematical entropy pair for stiff balance laws}

We now assuming now that a mathematical entropy 
$\Phi:\Omega\to\RR$ exists and satisfies the following two additional conditions: 

\

\noindent{{\bf  Condition 5.} There exists an entropy-flux $\Psi:\Omega\to\RR$ such that  
$D_U \Phi \,  A = D_U \Psi$ in $\Omega$.
So, all smooth solutions satisfy
$$
\eps\del_t \big( \Phi(U^\eps) \big) + \del_x \big( \Psi(U^\eps) \big)
= 
- D_U \Phi(U^\eps)\,  {R(U^\eps) \over \eps}
$$
and, consequently, the matrix $D^2_U \Phi \, A$ is symmetric in $\Omega$.
Moreover, the map $\Phi$ is convex, i.e.~the $N\times N$ 
matrix $D^2_U \Phi$ is positive definite on $\Mcal$.

\

\noindent{{\bf  Condition 6.}  
The entropy is compatible withe the relaxation in the sense that   
$$
\aligned
&  D_U \Phi \,  R \ge 0 \qquad  \text{ in } \Omega, 
\\
& D_U \Phi(U) = \nu(U) Q \in \RN, \qquad \nu(U) \in \Rn.  
\endaligned
$$
   
\

Next, we return to the effective equations 
$\del_t u= \del_x \Dcal$ sand 
$\Dcal := -Q  A(\Ecal(u))\,  U_1 $
and, multiplying it by the Hessian of the entropy, we see that  
$U_1\in\RR^N$ is characterized by 
$$
\aligned
\Lcal(u) U_1 &= - \big( D^2_U \Phi \big) (\Ecal(u)) \del_x \big( F(\Ecal(u)) \big),
\\
Q U_1 &= 0, 
\endaligned
$$
with $\Lcal(u) = D^2_U \Phi(\Ecal(u))B(\Ecal(u))$. 

Denoting by $\Lcal(u)^{-1}$ the generalized inverse with constraint and setting
$S(u) := Q A(\Ecal(u))$,  
we obtain  
$$
\Dcal = S \Lcal^{-1} \big( D^2_U \Phi \big) (\Ecal) \del_x\big( F(\Ecal) \big).
$$
Finally, one can check that, with 
$v := \del_x \left(  D_u \Phi(\Ecal) \right)^T $,  
$$
\big( D^2_U \Phi \big) (\Ecal) \del_x \big( F(\Ecal) \big) = S^T  v.
$$

\begin{theorem}[Entropy structure of the effective system]
When a mathematical entropy is available, the effective equations take the form 
$$
\del_t u = \del_x \Big( L(u) \, \del_x \big( D_u \Phi(\Ecal(u)) \big)^T \Big), 
$$
with 
$$
\aligned
L(u) & := S(u)\Lcal(u)^{-1}S(u)^T, 
&
S(u) := Q A(\Ecal(u)), 
\\
\Lcal(u) &:= \big( D^2_U\Phi \big) (\Ecal(u)) B (\Ecal(u)),
\endaligned 
$$
where, for all $b$ satisfying $Qb =0$,  the unique solution to
$\Lcal(u) V =b$, $QV =0$
is denoted by $\Lcal(u)^{-1} b$ (generalized inverse).  
\end{theorem}

This result can be formulated in the so-called entropy variable 
$\big( D_u \Phi(\Ecal(u)) \big)^T$. 
Furthermore, a dissipation property follows from our assumptions and, specifically, from the entropy and equilibrium properties (see $R(\Ecal(u))=0$),  we obtain 
$$
\aligned
D_U \Phi R & \ge 0 \qquad  \text { in } \Omega,
\\
\left(D_U \Phi R \right)|_{ U=\Ecal(u)}& =0 \qquad  \text {in } \omega.
\endaligned
$$
Thus, the matrix $D^2_U\Big( D_U \Phi R \Big)|_{ U=\Ecal(u)}$ is non-negative definite. 
 It follows that 
$$ 
D^2_U\Big(D_U \Phi R \Big) 
 = D^2_U\Phi B + \left( D^2_U\Phi B \right)^T \qquad \text{ when } U=\Ecal(u),
$$
so that 
$D^2_U\Phi \, B |_{ U=\Ecal(u)} \geq 0$ in $\omega$.

For the {equilibrium entropy} $\Phi(\Ecal(u))$, the associated 
(entropy) flux $u \mapsto \Psi(\Ecal(u))$ is constant on the equilibrium manifold $\omega$. 
For the map $\Psi(\Ecal)$, we have 
$$
\aligned
D_u \big( \Psi(\Ecal) \big) & = D_U \Psi(\Ecal) D_u \Ecal
 = D_U \Phi(\Ecal) A(\Ecal) D_u \Ecal.
\endaligned
$$
Observe that $\big( D_U \Phi \big) (\Ecal) = D_u \big( \Phi(\Ecal)\big) Q$, so that  
$$
\aligned
D_u \Big( \Psi(\Ecal(u)) \Big) &= D_u \Phi(\Ecal(u)) Q  A(\Ecal(u)) D_u
\Ecal(u)
\\
&=  D_u \Big( \Phi(\Ecal(u)) \Big) D_u QF(\Ecal(u)). 
\endaligned
$$
Since $QF(\Ecal) = 0$, then $D_u QF(\Ecal) =0$ and the proof is completed.

Therefore, $D_u \big( \Psi(\Ecal(u)) \big)=0$ for all $u \in \omega$. From the expansion
$U^\eps = \Ecal(u) + \eps U_1 + ...$, 
where $U_1$ is given by the first-order corrector problem, we deduce 
$$
\Psi(U^\eps)=\Psi(\Ecal(u)) +\eps \, D_U\Psi(\Ecal(u)) \, U_1 + \Ocal(\eps^2),
$$
and then 
$\del_x \Psi(U^\eps) = \eps \, \del_x D_U\Psi(\Ecal(u)) \, U_1 + \Ocal(\eps^2)$. 
Similarly, for the relaxation source,  we have 
$$
D_U \Phi(U^\eps)R(U^\eps) = \eps^2D^2_U \Phi(\Ecal(u)) D_U
R(\Ecal(u)) U_1 + \Ocal(\eps^3).
$$
We thus get 
$$
\aligned
& \del_t \big( \Phi(\Ecal(u)) \big)+ \del_x \Big( D_U \Psi(\Ecal(u)) \, U_1\Big) 
\\
& = - U_1^T \,  \left( D^2_U \Phi(\Ecal(u)) B(\Ecal(u)) \right) U_1. 
\endaligned
$$ 
At this juncture, recall that 
$X\,  \, \big( D^2_U \Phi \big)(\Ecal) B(\Ecal) \, X \ge 0$ for $X \in \RR^N$.  

\begin{proposition}[Monotonicity of the entropy]
The entropy is non-increasing, i.e. 
$$
\del_t \big( \Phi(\Ecal(u)) \big) + \del_x \left( D_U \Psi(\Ecal(u)) \, U_1\right) \leq 0 
$$
and  
$$
\del_t \big( \Phi(\Ecal(u)) \big) = \del_x \Big( \big( D_u \big(\Phi(\Ecal(u) \big) \big) 
L(u) \, \del_x \big( D_u \big(\Phi(\Ecal(u) \big) \big)^T \Big).  
$$
\end{proposition}


\subsection{Effective models} 

\subsubsection*{Effective model for stiff friction}

We now analyze the diffusive regime for the Euler equations with friction. According to the general theory, the equilibria satisfy 
$\del_t\rho = -\del_x\Big( Q A(\Ecal(u)) \, U_1 \Big)$ 
with 
$$
D_U F(\Ecal(u)) = \left(\begin{array}{cc}
0 & 1 \\ p'(\rho) & 0
\end{array}\right).
$$
Here, $U_1$ is the unique solution to
$B (\Ecal(u)) U_1 = -\del_x \big( F(\Ecal(u)) \big)$ and 
$QU_1 = 0$
with 
$$
B (\Ecal(u)) = \left(\begin{array}{cc}
0 & 0 \\ 0 & 1
\end{array}\right), 
\quad\quad
\del_x \Big( F(\Ecal(u)) \Big) = \left(\begin{array}{c}
0 \\ \del_x \big( p(\rho)\big) 
\end{array}\right). 
$$
The effective diffusion equation for the Euler equations with friction thus read :  
\be
\label{Euler2}
\del_t\rho = \del^2_x \big( p(\rho) \big), 
\ee
which is a nonlinear parabolic equation (away from vacuum) since $p'(\rho) >0$. 
Near the vacuum, this equation is often degenerate since  $p'(\rho)$ typically vanishes at $\rho=0$. 
For instance for polytropic gases $p(\rho) = \kappa \rho^\gamma$ with 
$\kappa>0$ and $\gamma \in (1,\gamma)$ we get 
\be
\label{Euler3}
\del_t\rho = \kappa \gamma \, \del_x \big( \rho^{\gamma-1} \del \rho \big). 
\ee

Defining the internal energy $e(\rho)>0$ by 
$e'(\rho)= p(\rho) / \rho^2$
we see that, for all smooth solutions to \eqref{Euler1}, 
\be
\label{Euler3-2} 
\eps \, \del_t \Big(\rho\frac{v^2}{2}+\rho e(\rho)\Big)
+
\del_x\Big(\rho\frac{v^3}{2}+ (\rho e(\rho) + p(\rho)) v \Big)
=-\frac{\rho v^2}{\eps}, 
\ee
so that $\Phi(U)=\rho\frac{v^2}{2}+\rho e(\rho)$
 is a convex entropy and is compatible with the relaxation. 
All the conditions of the general framework are therefore satisfied.


\subsubsection*{Effective model for stiff radiative transfer}
 
This system is      compatible with our late-time/stiff relaxation framework with now 
$$
U=\left(\begin{array}{c}
e \\ f \\ \tau
\end{array}\right),
\qquad
F(U)=\left(\begin{array}{c}
f \\ \chi(\frac{f}{e})e \\ 0
\end{array}\right), 
\quad\quad
R(U)=\left(\begin{array}{c}
e-\tau^4 \\ f \\ \tau^4-e
\end{array}\right).
$$
The equilibria read $u=\tau+\tau^4$ and 
$$
\Ecal(u)=\left(\begin{array}{c}
\tau^4 \\ 0 \\ \tau
\end{array}\right),
\qquad \quad
Q := (1~0~1)
$$ 
and we have $QF(\Ecal(u))=0$.  

We determine the diffusive regime for the $M1$ model from   
$$
\big( D_U F \big) (\Ecal(u)) = 
\left(\begin{array}{ccc}
0 & 1 & 0 \\ 
\chi(0) & \chi'(0) & 0 \\
0 & 0 & 0
\end{array}\right)
=
\left(\begin{array}{ccc}
0 & 1 & 0 \\ 
\frac{1}{3} & 0 & 0 \\
0 & 0 & 0
\end{array}\right), 
$$
where $U_1$ is the solution to 
$$
\begin{aligned}
\left(\begin{array}{ccc}
1 & 0 & -4\tau^3 \\ 
0 & 1 & 0 \\
-1 & 0 & 4\tau^3
\end{array}\right)U_1 & =
\left(\begin{array}{c}
0 \\ 
\del_x\big( \tau^4/3 \big) \\
0 
\end{array}\right),\\
 (1~0~1) U_1 & =0. 
\end{aligned}
$$
Therefore, we have  
$U_1 = \left(\begin{array}{c}
0 \\ 
\frac{4}{3} \tau^3\del_x\tau \\
0
\end{array}\right)$ 
and the effective diffusion equation reads  
\be
\label{M12} 
\del_t(\tau+\tau^4) = \del_x\left(
\frac{4}{3}\tau^3\del_x \tau \right), 
\ee
which      admits an entropy.  


\subsubsection*{Effective model for stiff friction and stiff radiative transfert}
  
Here, we have 
$$
D_U F(\Ecal(u))=\begin{pmatrix}0&1&0&0\\p'(\rho)&0&0&0\\0&0&0&1\\0&0&\frac{1}{3}&0\end{pmatrix},
\qquad
U_1=\begin{pmatrix}0\\ \frac{1}{\kappa}\Bigl(-\del_x p(\rho)-\frac{1}{3}\del_x e\Bigr)\\0
\\-\frac{1}{3\sigma}\del_x e\end{pmatrix}, 
$$ 
and the effective diffusion system for the coupled Euler/$M1$ model reads  
\be
\label{EM2} 
\begin{aligned}
&\del_t \rho-\frac{1}{\kappa}\del_x^2 p(\rho)-\frac{1}{3\kappa}\del_x^2 e=0,\\
&\del_t e-\frac{1}{3\sigma}\del_x^2 e=0.
\end{aligned}
\ee
The second equation is a heat equation, and its solution appears as a source-term in the first one. 

 
\subsubsection*{Effective model for stiff nonlinear friction}

Our framework      encompass handle certain nonlinear diffusion regime under the scaling 
$$
\eps \, \del_t U+\del_x F(U) = -\frac{ R(U)}{\eps^q}.
$$
The parameter $q \ge 1$ introduces a new scale and is necessary when the relaxation is nonlinear. 
We assume that 
$$
R\big(\Ecal(u)+\eps \, U\big) = \eps^q R\Big(\Ecal(u)+M(\eps) \, U\Big),
\qquad U\in\Omega, \quad u\in\omega 
$$
for some matrix $M(\eps)$.  In that regime, the effective equations are now nonlinear parabolic.   
 
Our final example requires this more general theory and reads 
\be
\label{EN1} 
\begin{aligned}
& \eps \del_t h +\del_x \big( hv \big) 
=0,\\
& \eps \del_t \big( hv \big) 
+\del_x \Big( h \, v^2 + p(h) \Big) = -\frac{\kappa^2(h)}{\eps^2} \, g \, hv|hv|,
\end{aligned}
\ee
where $h$ is the fluid height and $v$ the fluid velocity $v$. The pressure reads $p(h)=g\, h^2/2$ while $g>0$ is the gravity constant. The friction  $\kappa:\RR^+ \to\RR^+$ is a positive function, and for instance one can take  
$\kappa(h)=\frac{\kappa_0}{h}$ with $\kappa_0>0$.

The nonlinear version of the late-time/stiff relaxation framework applies by introducing 
$$
U=\left(\begin{array}{c}
h \\ hv
\end{array}\right),
\quad
F(U)=\left(\begin{array}{c}
hv \\ hv^2+p(h)
\end{array}\right), 
\quad
R(U)=\left(\begin{array}{c}
0 \\ \kappa^2(h)ghv|hv|
\end{array}\right).
$$
The equilibria $u=h$ are associated with  
$$
\Ecal(u)=\left(\begin{array}{c}
h \\ 0
\end{array}\right), 
\qquad \qquad 
Q=(1~0). 
$$ 
The relaxation is nonlinear and  
$$
R(\Ecal(u)+\eps U) = \eps^2 R\big( \Ecal(U) + M(\eps)U \big),
$$
with 
$$
M(\eps) := \left(\begin{array}{cc}
\eps & 0 \\ 0 & 1
\end{array}\right).
$$ 
in turn, we obtain a nonlinear effective equation for the Euler equations with nonlinear friction, i.e. 
\be
\label{EN2} 
\del_t h = \del_x\left( {\sqrt{h} \over \kappa(h)} \, 
\frac{\del_x h}{\sqrt{|\del_x h|}}\right),
\ee
which is a parabolic and fully nonlinear. 

Introducing the internal energy $e(h) :=gh/2$, we see that all smooth solutions to \eqref{EN1} 
satisfy the entropy inequality 
\be
\label{EN2-2} 
\eps\del_t\left( h\frac{v^2}{2}+g\frac{h^2}{2} \right)
+
\del_x\left(
h\frac{v^2}{2}+gh^2
\right)v=-\frac{\kappa^2(h)}{\eps^2}ghv^2|hv|.
\ee
The entropy 
$\Phi(U) := h\frac{v^2}{2}+g\frac{h^2}{2}$ 
satisfies the compatibility properties for the nonlinear late-time/stiff relaxation theory, with 
$$
R(\Ecal(u)+M(0)\bar{U}_1) = \left(\begin{array}{c}
0 \\ \del_x p(h)
\end{array}\right),
$$
where $\bar{U}_1=(0~\beta)\, $.   We obtain $R(\Ecal(u)+M(0)\bar{U}_1)=c(u)\bar{U}_1$ with
$$
c(u) = g\kappa(h)\sqrt{h|\del_x h|}\ge 0.
$$


\subsection{A class of asymptotic-preserving finite volume method} 

\subsubsection*{The general strategy}

We now will design a class of finite volume schemes which are consistent with the asymptotic regime $\eps \to 0$
and allow us to recover the effective diffusion equation
(independently of the  mesh-size) for the limiting solutions. Hence, we develop here a rather
 general framework adapted to the hyperbolic-to-parabolic relaxation regime.

\begin{itemize}

\item[Step 1.]  \qquad We rely on a arbitrary finite volume scheme for  the homogeneous system
$$ 
\del_t U+\del_x F(U) =0, 
$$
as described below.  

\item[Step 2.]  \qquad  Next, we modify this scheme and include a matrix-valued free parameter in order to consistently 
approximate the non-homogeneous system  (for any $\gamma>0$) 
$$
\del_t U+\del_x F(U) = -\gamma \, R(U). 
$$

\item[Step 3.]  \qquad  By performing an asymptotic analysis of this scheme after replacing the discretization parameter
 $\Delta t$ by $\eps \Delta t$, and 
$\gamma$ by $1/\eps$, 
we then determine the free parameters and ensure the desired asymptotic-preserving property.
\end{itemize}

For definitness, the so-called HLL discretization of the homogeneous system  
(Harten, Lax, and van Leer \cite{HLL}) are now discussed. We present the solver based on 
a single intermediate state and on a uniform mesh with cells of length $\Delta x$, that is,  
$$
[x_{i-\demi},x_{i+\demi}], 
\qquad 
x_{i+\demi}=x_i+\frac{\Delta x}{2} 
$$
for all $i=\ldots, -1,0,1, \ldots$.  
The time discretization is based on some $\Delta t$ restricted by the CFL condition~\cite{CFL} with  
$
t^{m+1}=t^m+\Delta t. 
$

Given any initial data (lying in $\Omega$): 
$$
U^0(x) = \frac{1}{\Delta x}
\int_{x_{i-\demi}}^{x_{i+\demi}} U(x,0)dx,
\qquad x\in[x_{i-\demi},x_{i+\demi}).
$$
we design approximations that are piecewise constant at each $t^m$, that is,  
$$
U^m(x)=U_i^m, 
\qquad x\in[x_{i-\demi},x_{i+\demi}), \quad i\in\ZZ. 
$$

At each cell interface we use the {approximate Riemann solver} 
$$
\Ubar_\Rcal(\frac{x}{t};U_L,U_R) =\left\{
\begin{aligned}
& U_L, \quad \frac{x}{t}<-b,\\
& \Ubar^\star, \quad -b<\frac{x}{t}<b,\\
& U_R, \quad \frac{x}{t}>b,
\end{aligned}\right.
$$
where $b>0$ is (sufficiently) large.
The ``numerical cone'' (and numerical diffusion) is determined by some  $b>0$ and, for simplicity in the presentation, we assume
 a single constant $b$. More generally,  one can introduce distinct speeds $b^-_{i+\demi}<b^+_{i+\demi}$ at each interface. 

We introduce the intermediate state 
$$
\Ubar^\star = \frac{1}{2}(U_L+U_R)-\frac{1}{2b}
\big( F(U_R) - F(U_L) \big) 
$$
and, under the CFL condition 
$b\frac{\Delta t}{\Delta x}\le 1/2$, 
the underlying Riemann solutions are non-interacting. Our global approximations 
$$
\Ubar^m_{\Delta x}(x,t^m+t), \qquad t\in[0,\Delta t), \quad x \in \RR. 
$$
are defined as follows. 

At the time $t^{m+1}$, we set  
$$
\Ubar_i^{m+1} = \frac{1}{\Delta x}
\int_{x_{i-\demi}}^{x_{i+\demi}} \Ubar^m_{\Delta x}(x,t^m+\Delta t) dx
$$
and, recalling 
$
\Ubar^\star_{i+\demi} = \frac{1}{2}(U_i^m+U_{i+1}^m)-\frac{1}{2b}
(F(U_{i+1}^m)-F(U_i^m)),
$
and integrating out the expression given by the Riemann solutions, 
we arrive at the scheme adapted to our homogeneous system 
$$
\Ubar_i^{m+1} = U_i^m - \frac{\Delta t}{\Delta x}
\Big( F^{HLL}_{i+\demi} - F^{HLL}_{i-\demi} \Big),
$$
where
$$
F^{HLL}_{i+\demi} = \frac{1}{2} \Big( F(U_i^m)+F(U_{i+1}^m) \Big)
-\frac{b}{2}(U_{i+1}^m-U_i^m). 
$$
More generally one can include here two speeds $b^-_{i+\demi}<b^+_{i+\demi}$.

This scheme enjoys an invariant domain property, as follows. The intermediate states $\Ubar^\star_{i+\demi}$ can be written in the form of a convex combination 
$$
\Ubar^\star_{i+\demi}= 
\frac{1}{2}\left(U_i^m+\frac{1}{b}F(U_i^m)\right)
+ 
\frac{1}{2}\left(U_{i+1}^m-\frac{1}{b}F(U_{i+1}^m)\right) \in \Omega,
$$
provided $b$ is large enough. An alternative decomposition is 
$$
\Ubar^\star_{i+\demi}
= 
\frac{1}{2}\Big( I + \frac{1}{b} \Abar(U_i^m, U^m_{i+1}) \Big) \, U_i^m
+ 
\frac{1}{2}\Big( I - \frac{1}{b} \Abar(U_i^m, U^m_{i+1}) \Big) \, U_{i+1}^m, 
$$
where $\Abar$ is an  ``average'' of $D_U F$.   
By induction, we conclude that  $\Ubar_i^m$ in $\Omega$ for all $m,i$.


\subsubsection*{Handling the stiff relaxation} 

Consider the {modified Riemann solver:} 
$$
U_\Rcal(\frac{x}{t};U_L,U_R) =\left\{
\begin{aligned}
& U_L, \quad \frac{x}{t}<-b,\\
& U^{\star L}, \quad -b<\frac{x}{t}<0,\\
& U^{\star R}, \quad 0<\frac{x}{t}<b,\\
& U_R, \quad \frac{x}{t}>b,
\end{aligned}\right.
$$
with, at the interface,  
$$
\begin{aligned}
& U^{\star L} = \alphab\Ubar^\star +  (I-\alphab)\big( U_L - \bar{R}(U_L) \big),
\\
& U^{\star R} = \alphab\Ubar^\star + 
    (I-\alphab) \big( U_R - \bar{R}(U_R) \big).
\end{aligned}
$$
We have introduced an arbitrary $N\times N$-matrix and an $N$-vector by 
$$
\alphab = \left( I +\frac{\gamma\Delta
  x}{2b}(I+\underline{\sigma})
\right)^{-1},
\qquad  
\bar{R}(U) = (I+\underline{\sigma})^{-1} R(U). 
$$
The term $\underline{\sigma}$ is a parameter matrix and we require that  all inverse matrices are well-defined and, importanty,
 the correct asymptotic regime arises at the discrete level (see below).

At each  $x_{i+\demi}$, we use the Riemann solver 
$U_{\Rcal}(\frac{x-x_{i+\demi}}{t-t^m}; U_i^m,U_{i+1}^m)$ 
and superimpose non-interacting Riemann solutions 
$$
U^m_{\Delta x}(x,t^m+t), \qquad t\in[0,\Delta t), \quad x \in \RR.
$$ 
The approximation at the time $t^{m+1}$  reads  
$U_i^{m+1} = 
\int_{x_{i-\demi}}^{x_{i+\demi}} U^m_{\Delta x}(x,t^m+\Delta t)dx$. 
By integration of the Riemann solutions, we arrive at the following discrete form of the balance law 
\be
\label{904} 
\begin{aligned}
& \frac{1}{\Delta t}(U_i^{m+1}-U_i^m) +\frac{1}{\Delta x}
\Big(\alphab_{i+\demi} F^{HLL}_{i+\demi}
-\alphab_{i-\demi} F^{HLL}_{i-\demi} \Big)
\\
& =
\frac{1}{\Delta x}
  (\alphab_{i+\demi}-\alphab_{i-\demi}) F(U_i^m)
-\frac{b}{\Delta x} (I-\alphab_{i-\demi})
  \bar{R}_{i-\demi}(U_i^m)
\\
& \quad - 
\frac{b}{\Delta x} (I-\alphab_{i+\demi})
  \bar{R}_{i+\demi}(U_i^m). 
\end{aligned}
\ee

The source can rewritten as 
$$
\aligned
\frac{b}{\Delta x}(I-\alphab_{i+\demi})
  \bar{R}_{i+\demi}(U_i^m) 
& = 
\frac{b}{\Delta x}\alphab_{i+\demi}
 (\alphab_{i+\demi}^{-1}-I)
  \bar{R}_{i+\demi}(U_i^m)
\\
& = \frac{\gamma}{2}\alphab_{i+\demi}R(U_i^m) 
\endaligned
$$
and 
$$
\frac{b}{\Delta x}(I-\alphab_{i-\demi})
  \bar{R}_{i-\demi}(U_i^m) = 
\frac{\gamma}{2}\alphab_{i-\demi}R(U_i^m). 
$$
Our finite volume scheme for late-time/stiff-relaxation problems finally read  
\be
\label{914} 
\begin{aligned}
& \frac{1}{\Delta t}(U_i^{m+1}-U_i^m) + \frac{1}{\Delta x}
(\alphab_{i+\demi} F^{HLL}_{i+\demi}
-\alphab_{i-\demi} F^{HLL}_{i-\demi})
\\
& =
\frac{1}{\Delta x}
  (\alphab_{i+\demi}-\alphab_{i-\demi}) F(U_i^m)
-\frac{\gamma}{2}
(\alphab_{i+\demi}+\alphab_{i-\demi})R(U_i^m).
\end{aligned}
\ee

\begin{theorem}[A class of finite volume schemes for relation problems]
When 
$$
\underline{\sigma}_{i+\demi}-\underline{\sigma}_{i-\demi}=\Ocal(\Delta
x)
$$ 
and the matrix-valued map $\underline{\sigma}$ is smooth, 
the finite volume scheme above is {\rm consistent} with the hyperbolic system with relaxation and satisfies 
the following invariant domain property:  
provided all states 
$$
\begin{aligned}
&  U^{\star L}_{i+\demi} = \alphab_{i+\demi}\Ubar^\star_{i+\demi} + 
    (I-\alphab_{i+\demi})(U_i^m-\bar{R}(U_i^m)),\\
& U^{\star R}_{i+\demi} = \alphab_{i+\demi}\Ubar^\star_{i+\demi} + 
    (I-\alphab_{i+\demi})(U_{i+1}^m-\bar{R}(U_{i+1}^m))
\end{aligned}
$$
belong to $\Omega$, then all of the states $U_i^m$ belong to $\Omega$.
\end{theorem}


\subsection{Effective equation for the discrete asymptotics}

We replace $\Delta t$ by $\Delta t/\eps$ and $\gamma$ by $1 / \eps$ and consider the expression 
$$
\begin{aligned}
& \frac{\eps}{\Delta t}(U_i^{m+1}-U_i^m) +\frac{1}{\Delta x}
(\alphab_{i+\demi} F^{HLL}_{i+\demi}
-\alphab_{i-\demi} F^{HLL}_{i-\demi})
\\
&
= \frac{1}{\Delta x}
  (\alphab_{i+\demi}-\alphab_{i-\demi}) F(U_i^m)
-\frac{1}{2\eps}
(\alphab_{i+\demi}+\alphab_{i-\demi})R(U_i^m),
\end{aligned}
$$
in which 
$$
\alphab_{i+\demi}=\left(
I+\frac{\Delta x}{2\eps b}(I+\underline{\sigma}_{i+\demi})
\right)^{-1}.
$$
We expand near an equilibrium state 
$U_i^m = \Ecal(u_i^m) + \eps(U_1)_i^m + \Ocal(\eps^2)$ 
and find 
$$
\aligned 
F^{HLL}_{i+\demi}  
& =
{1 \over 2} F\big( \Ecal(u_i^m) \big) + {1 \over 2} F\big( \Ecal(u_{i+1}^m) \big) 
-\frac{b}{2}\left( \Ecal(u_{i+1}^m) - \Ecal(u_i^m) \right)
+ \Ocal(\eps),
\\
\frac{1}{\eps}R(U_i^m) & = B(\Ecal(u_i^m)) (U_1)_i^m +
  \Ocal(\eps), 
\\
\alphab_{i+\demi}
& =
\frac{2b \eps}{\Delta x} \big(I+\underline{\sigma}_{i+\demi} \big)^{-1} + \Ocal(1).
\endaligned
$$

The first-order terms yield us 
$$
\begin{aligned}
& \frac{1}{\Delta t}(\Ecal(u_i^{m+1}) - \Ecal(u_i^m))
\\
& = - \frac{2b}{\Delta x^2}\left(
(I+\underline{\sigma}_{i+\demi})^{-1} F^{HLL}_{i+\demi}|_{\Ecal(u)}
-
(I+\underline{\sigma}_{i-\demi})^{-1} F^{HLL}_{i-\demi}|_{\Ecal(u)}
\right)
\\
& \quad + 
\frac{2b}{\Delta x^2} \left(
(I+\underline{\sigma}_{i+\demi})^{-1}
-
(I+\underline{\sigma}_{i-\demi})^{-1}
\right) F(\Ecal(u_i^m))
\\
&
\quad - \frac{b}{\Delta x} \left(
(I+\underline{\sigma}_{i+\demi})^{-1}
+
(I+\underline{\sigma}_{i-\demi})^{-1}
\right)B(\Ecal(u_i^m))(U_1)_i^m.
\end{aligned}
$$
Assuming here the existence of an $n\times n$ matrix
$\Mcal_{i+\demi}$ satisfying 
$$
Q \big( I+\underline{\sigma}_{i+\demi} \big)^{-1}=\frac{1}{b^2} \, \Mcal_{i+\demi} Q 
$$
and multiplying the equation above by $Q$, we get 
$$
\frac{1}{\Delta t}(u_i^{m+1}-u_i^m) =
-\frac{2}{b\Delta x^2}\left(
\Mcal_{i+\demi} Q F^{HLL}_{i+\demi}|_{\Ecal(u)} - 
\Mcal_{i-\demi} Q F^{HLL}_{i-\demi}|_{\Ecal(u)} \right),
$$
where 
$$
\begin{aligned}
Q F^{HLL}_{i+\demi}|_{\Ecal(u)}&=  {Q \over 2}  F(\Ecal(u_i^m)) 
+ {Q \over 2}  F(\Ecal(u_{i+1}^m))
-\frac{b}{2} Q \left(\Ecal(u_{i+1}^m)-\Ecal(u_i^m) \right)
\\
&= -\frac{b}{2} (u_{i+1}^m -u_i^m).
\end{aligned}
$$
The asymptotic system for the scheme thus reads 
\be
\label{944} 
\frac{1}{\Delta t}(u_i^{m+1}-u_i^m) =
\frac{1}{\Delta x^2}
\Big(
\Mcal_{i+\demi}(u_{i+1}^m -u_i^m) + 
\Mcal_{i-\demi}(u_{i-1}^m -u_i^m)
\Big).
\ee
Recall that for some matrix $\Mcal(u)$, the effective equation reads
$\del_t u = \del_x \left(\Mcal(u) \del_x u \right)$.

\begin{theorem} [Discrete  late-time asymptotic-preserving property] 
Assume that the matrix-valued coefficients satisfy the following conditions: 

\begin{itemize} 

\item[$\bullet$] The matrices $I+\underline{\sigma}_{i+\demi}$ 
and $\Big(1+\frac{\Delta x}{2\eps b} \Big)I+\underline{\sigma}_{i+\demi}$
are  invertible for all $\eps \in [0,1]$.

\item[$\bullet$] There exists a matrix $\Mcal_{i+\demi}$ satisfying the commutation condition 
$$
Q(I+\underline{\sigma}_{i+\demi})^{-1}=\frac{1}{b^2}\Mcal_{i+\demi} Q. 
$$

\item[$\bullet$] The discrete formulation of $\Mcal(u)$ at each interface $x_{i+\demi}$ satisfies 
$$
\Mcal_{i+\demi} = \Mcal(u) + \Ocal(\Delta x). 
$$
\end{itemize}
Then the effective system associated with the proposed finite volume scheme coincides with the effective system 
determined in the late-time/stiff relaxation framework. 
\end{theorem}

Finally, wee refer to \cite{BLT} for various numerical experiments demonstrating the relevance of the proposed scheme and its 
efficiency in order to compute late-time behaviors of solutions. Asymptotic solutions may have large gradients but are in fact regular. 
Note      that our CFL stability condition is based on the homogeneous hyperbolic system and therefore imposes a restriction 
on $\Delta t/\Delta x$ only.
In our test, for simplicity, the initial data were taken in the image of $Q$, while the reference solutions (needed for the purpose of comparison) 
were computed separately by solving the associated parabolic equations, of course under a (much more restrictive) restriction on $\Delta t/(\Delta x)^2$. 

The proposed theoretical framework for late-time/stiff relaxation problems thus led us 
to the development of a good strategy to design asymptotic-preserving schemes involving matrix-valued  parameter. 
The convergence analysis ($\eps \to 0$) and the numerical analysis ($\Delta x \to 0$) for the problems under consideration 
are important and challenging open problems. It would be      very interesting to apply our technique to plasma mixtures in a  multi-dimensional setting.  
 
Furthermore, high-order accurate Runge-Kutta methods have been recently developed for these stiff relaxation problems by Boscarino and Russo \cite{BR} and by Boscarino, LeFloch, and Russo \cite{BLR}.


\section{Geometry-preserving finite volume methods}
\label{sec:3} 

\subsection{Objective and background material}

On a smooth $(n+1)$-dimensional manifold $M$ refered to as a spacetime, we consider the class of {nonlinear conservation laws}
\be
\label{cons1}
d(\om(u)) = 0, \qquad u=u(x), \, x \in M.
\ee
For all $\ubar \in \RR$, $\omega=\omega(\ubar)$ is a smooth {field of $n$-forms}, refered to as the {flux field} of the conservation law under consideration.

Two examples are of particular interest. When $M = \RR_+ \times N$ and the $n$-manifold $N$ is 
endowed with a {Riemannian metric} $h$, \eqref{cons1} reads 
$$
\del_t u + \dive_h (b(u)) = 0, \qquad u=u(t,y), \ t \geq 0, \, y \in N, 
$$
where $\dive_h$ denotes the divergence operator for the metric $h$. The flux field is then considered as a {flux vector field} $b=b(\ubar)$ {on the $n$-manifold $N$} and is independent of the time variable. 

More generally, when $M$ is endowed with a {Lorentzian metric} $g$, \eqref{cons1} reads 
$$
\dive_g (a(u)) = 0, \qquad u=u(x), \, x \in M, 
$$
in which the flux $a=a(\ubar)$ is now a vector field {on $M$.}
In this Riemannian or Lorentzian settings, the theory of weak solutions on manifolds
was initiated by Ben-Artzi and LeFloch \cite{BL} and developed in \cite{ABL,ALO,PLF,LO}. 

In the present apaper, we discuss the novel approach in which the conservation law is written in the form \eqref{cons1}, that is, 
the flux $\omega=\omega(\ubar)$ is defined as a {field of differential forms of degree $n$.}
No geometric structure is assumed on $M$ and the sole flux field structure
is assumed. The equation \eqref{cons1} is a ``conservation law'' for the unknown quantity $u$, as follows from Stokes theorem 
for sufficiently smooth solutions $u$: the total flux
\be
\label{486}
\int_{\del \Ucal} \om(u) = 0, \qquad \Ucal \subset M,
\ee
vanishes for every smooth open subset $\Ucal$.
By relying on \eqref{cons1} rather than the equivalent expressions in the special cases 
of Riemannian or Lorentzian manifolds,
we develop a theory of entropy solutions 
which is technically and conceptually simpler and  provides a generalization of earlier works. From
 a numerical perspective, relying o \eqref{cons1} leads us to a geometry-consistent class of finite volume schemes, as we will now present it. 
So, our main objective i this presentation will be a generalization of the formulation and convergence of the 
finite volume method for general conservation law \eqref{cons1}. In turn, this will also 
establish the existence of a contracting semi-group of entropy solutions.  

 We will proceed as follows:

\bei 

\item First we will  formulate the initial and boundary problem for \eqref{cons1} by 
taking into account the nonlinearity and hyperbolicity 
of the equation. We need to impose
that the manifold satisfies a {global hyperbolicity condition,} which provides a global time-orientation
and allow us to distinguish between ``future'' and ``past'' directions in the time-evolution and we 
suppose that the manifold is foliated by compact slices.

\item Second, we introduce a geometry-consistent version of the finite volume method which provides a natural discretization of the conservation law \eqref{cons1}, which solely uses the $n$-volume form structure
associated with the flux field $\omega$.

\item Third, we derive stability estimates, especially certain discrete versions of the entropy inequalities. We obtain a uniform control of the
entropy dissipation measure, which, however, is not sufficient by itself to establish the compactness of the sequence of solutions.
Yet, these stability estimates imply that the sequence of approximate solutions generated
by the finite volume scheme converges to an {entropy measure-valued solution} in the sense of DiPerna.

\item Fourth, to conclude we rely on DiPerna's uniqueness theorem \cite{DiPerna} and 
establish the existence of entropy solutions to the corresponding initial value problem.

\eei

In the course of our analysis, we  will derive the following {contraction property}: 
for any entropy solutions $u, v$ and any hypersurfaces $H, H'$ such that $H'$ lies in the future of $H$, one has 
\be
\label{489}
\int_{H'} \Omegabf(u_{H'}, v_{H'}) \leq \int_{H} \Omegabf(u_{H}, v_{H}).
\ee  
Here, for all reals $\ubar, \vbar$, the $n$-form field
$\Omegabf(\ubar, \vbar)$ is determined from the flux field $\omega(\ubar)$
and is a generalization (to the spacetime setting) of the notion (introduced in \cite{Kruzkov}) 
of Kruzkov entropy $|\ubar - \vbar|$.

DiPerna's measure-valued solutions were first used to establish the convergence of schemes by Szepessy \cite{Szepessy},
 Coquel and LeFloch \cite{CL1,CL2,CL3}, and Cockburn, Coquel, and LeFloch \cite{CCL00,CCL2}.  
Further hyperbolic models including      a coupling with elliptic equations and many 
applications were 
investigated by Kr\"oner \cite{Kroener}, and Eymard, 
Gallouet, and Herbin \cite{EGH}. For higher-order schemes, see Kr\"oner, Noelle, and Rokyta \cite{KNR}. 
See also Westdickenberg and Noelle \cite{WN}.   

\subsection{Entropy solutions to conservation laws posed on a spacetime}
\label{deux}
 
We assume that $M$ is an oriented, compact, differentiable $(n+1)$-manifold
with boundary.  Given an $(n+1)$-form $\alpha$, its {modulus} is defined as the $(n+1)$-form
$|\alpha| : = |\overline{\al}|\, dx^0 \wedge \cdots \wedge dx^n$, 
where $\alpha = \overline{\al} \,dx^1 \wedge \cdots \wedge dx^n$ is written in an oriented frame
determined in coordinates $x=(x^\alpha)=(x^0, \ldots, x^n)$.
If $H$ is a hypersurface, we denote by $i=i_H : H \to \MM$ the canonical injection map, and
by $i^*=i_H^*$ is the pull-back operator acting on differential forms defined on $M$.

We introduce the following notion:
\bei

\item A {\rm flux field} $\omega$ on the $(n+1)$-manifold $M$ is a parametrized family $\omega(\ubar) \in \Lambda^n(M)$ of smooth fields
of differential forms of degree $n$,
that depends smoothly upon the real parameter $\ubar$.

\item The {\rm conservation law} associated with a flux field $\om$ and with unknown $ u:  M \to \RR$ is
\be
\label{LR.1}
d\big(\om(u)\big)=0,
\ee
where $d$ is the exterior derivative operator and, therefore, $d\big(\om(u)\big)$ is a field
of differential forms of degree $(n+1)$.

\item A flux field $\omega$ is said to {\rm grow at most linearly}
if for every $1$-form $\rho$ on $M$
\be
\label{LR.2}
\sup_{\ubar \in \RR} \int_M \left| \rho \wedge \del_u \om(\ubar) \right| < +\infty.
\ee
\eei

In local coordinates $x=(x^\al)$ we write (for all $\ubar \in \RR$)
$\om(\ubar) = \om^\al(\ubar) \,  (\widehat{dx})_\al$
and $(\widehat{dx})_\al  := dx^0 \wedge \ldots \wedge dx^{\al-1} \wedge dx^{\al+1} \wedge \ldots \wedge dx^n$. 
Here, the coefficients $\omega^\alpha = \om^\al(\ub)$ are smooth.
The operator $d$ acts on differential forms and that, given a $p$-form $\rho$ and a $p'$-form $\rho'$, one has $d(d\rho)=0$ and
$d(\rho \wedge \rho')= d\rho \wedge \rho' +(-1)^p \rho \wedge d\rho'$. The equation \eqref{LR.1} makes sense for unknowns that are  
Lipschitz continuous. However, solutions to nonlinear hyperbolic equations
need not be continuous and we need to recast \eqref{LR.1} in a weak form.

Given a {smooth} solution $u$ of \eqref{LR.1} we apply Stokes theorem on any open subset $\Ucal$
(compactly included in $M$ and with smooth boundary $\del \Ucal$) and find 
\be
\label{key64}
0 = \int_\Ucal d(\om(u)) = \int_{\del \Ucal} i^*(\om(u)).
\ee
Similarly, given any smooth function $\psi: M \to \RR$ we write
$d(\psi \, \om(u)) = d\psi \wedge \om(u) + \psi \, d(\om(u))$, 
where $d\psi$ is a $1$-form field. Provided $u$ satisfies \eqref{LR.1}, we deduce that 
$$
\int_M d(\psi \, \om(u)) = \int_M d\psi \wedge \om(u)
$$
and, by Stokes theorem,
\be
\label{LR.0}
\int_M d\psi\wedge \om(u)=\int_{\del \MM}i^*(\psi\om(u)).
\ee
A suitable orientation of the boundary $\del M$ is required for this formula to hold.

\begin{definition}[Weak solutions on a spacetime]
Given a flux field (with at most linear growth) $\omega$,
a function $u \in L^1(M)$ is a {\rm weak solution} to \eqref{LR.1} on the spacetime $M$
if
$\int_M d\psi\wedge \om(u) = 0$
for every $\psi : M \to \RR$ that is compactly supported in the interior $\mathring M$.
\end{definition}

Observe that the function $u$ is integrable and $\om(\ub)$ has at most linear growth in $\ub$,
so that the $(n+1)$-form $d\psi\wedge \om(u)$ is integrable on the compact manifold $M$.  

 
\begin{definition}
\label{key63}
A (smooth) field of $n$-forms $\Om=\Om(\ubar)$ 
is a {\rm (convex) entropy flux field} for \eqref{LR.1} if there exists
a (convex) function $U: \RR \to \RR$ such that
$$
\Om(\ubar) = \int_0^\ubar \del_u U (\vb) \, \del_u \om(\vb) \, d\vb, \qquad \ubar \in \RR.
$$
It is    {\rm admissible} if, moreover, $\sup | \del_u U | < \infty$.
\end{definition}

If we choose the function $U(\ub, \vb):=|\ub - \vb|$, where $\vb$ is a real parameter, the entropy flux field reads
\be
\label{KRZ}
\Omegabf(\ub, \vb) := \sgn(\ub - \vb) \, ( \om(\ub) - \om(\vb)). 
\ee
This is a generalization to spacetimes of the so-called Kruzkov's entropy pairs.
 
Next, given any smooth solution $u$ to \eqref{LR.1}, we multiply \eqref{LR.1} by $\del_u U(u)$
and obtain the conservation law
$$
d( \Om(u) ) -(d\Om)(u) + \del_u U(u) (d\om)(u) = 0.
$$
For discontinuous solutions, we impose the entropy inequalities
\be
\label{LR.1i}
d( \Om(u)) - (d\Om)(u) + \del_u U (u) (d\om)(u) \leq 0
\ee
 in the sense of distributions for all admissible entropy pair $(U,\Om)$.
This is justified, for instance, via the vanishing viscosity method, i.e.~by searching for weak solutions
 realizable as limits of smooth solutions to a parabolic
regularization.

It remains to prescribe initial and boundary conditions. We emphasize that, without further assumption
on the flux field (to be imposed shortly below), points along the boundary
$\del M$ can not be distinguished and it is natural to prescribe the trace of the solution
along the {whole} of the boundary $\del M$. This is possible provided the boundary data, $u_B: \del M \to \RR$,
is assumed by the solution in a suitably {weak sense}. Following Dubois and LeFloch \cite{DL}, 
we use the notation 
\be
u \big|_{\del M} \in \Ecal_{U,\Om}( u_B)
\label{LR.2i}
\ee
for all convex entropy pair $(U,\Om)$, where for all reals $\ub$ 
$$
\Ecal_{U,\Om} (\ub):= \Big\{ \vb \in \RR \, \, \big| \, \, E(\ub, \vb):= \Om(\ub) + \del_uU(\ub)(\om(\vb) - \om(\ub))
\leq  \Om(\vb)  \Big \}.
$$

\begin{definition}[Entropy solutions on a spacetime with boundary]
\label{key92}
Let $\om=\om(\ub)$ be a flux field (with at most linear growth) and
let $u_B \in L^1(\del M)$ be a boundary function.
A function $u \in L^1(M)$ is an {\rm entropy solution} to the boundary value problem $\eqref{LR.1}$
and $\eqref{LR.2i}$
if there exists a bounded and measurable field of $n$-forms $\gamma \in L^1\Lambda^n(\del M)$
such that, 
for every admissible convex entropy pair $(U,\Om)$ and every smooth function $\psi: M \to \RR_+$, 
$$
\aligned
&\int_M \Big(  d \psi \wedge \Om(u)
+ \psi \, (d  \Om) (u) - \psi \, \del_u U(u) (d \om) (u) \Big)
\\
& + \int_{\del \MM} \psi_{|\del M} \, \big(  i^*\Om(u_B) + \del_u U(u_B)  \big(\gamma - i^*\om(u_B) \big) \big)
\,  \geq 0. 
\endaligned
$$
\end{definition}

This definition makes sense since each of the terms
$d \psi \wedge \Om(u)$,  $(d  \Om) (u)$,  $(d \om) (u)$ belong to $L^1(M)$. 
Following DiPerna \cite{DiPerna}, we cn also consider solutions that are no longer functions but
{Young measures,} i.e,
weakly measurable maps $\nu: M \to \text{Prob}(\RR)$ taking values
within is the set of probability measures $\text{Prob}(\RR)$.

\begin{definition}
\label{LR.4}
Let $\om=\om(\ub)$ be a flux field with at most linear growth and let 
$u_B \in L^\infty(\del M)$ be a boundary function.
A compactly supported Young measure $\nu: M \to \text{Prob}(\RR)$ is
an {\rm entropy measure-valued solution}
to the boundary value problem $(\ref{LR.1}), (\ref{LR.2i})$
if there exists a bounded and measurable field of $n$-forms $\gamma \in L^\infty\Lambda^n(\del M)$
such that, for all convex entropy pair $(U,\Om)$ and all smooth functions $\psi \geq 0$, 
$$
\aligned
& \int_M  \Big\la \nu,  d \psi \wedge \Om(\cdot)  +  \psi \, \big(d  (\Om (\cdot)) - \del_u U(\cdot) (d \om) (\cdot)\big) \Big\ra
\\
& + \int_{\del M} \psi_{|\del M} \, \Big\la \nu, \Big( i^*\Om(u_B) + \del_u U(u_B)  \big(\gamma - i^*\om(u_B)\big)\Big) \Big\ra
\,  \geq 0. 
\endaligned
$$
\end{definition}


\subsection{Global hyperbolicity and geometric compatibility}
 
The manifold $M$ is now assumed to be foliated by hypersurfaces, say
\be
\label{foli}
M = \bigcup_{0 \leq t \leq T} H_t,
\ee
where each slice has the topology of a (smooth)
$n$-manifold $N$ with boundary. Topologically we have $M \simeq [0,T] \times N$, and 
\be
\label{foli2}
\aligned
\del M & = H_0 \cup H_T \cup B,
\\
B = (0,T) \times N & :=  \bigcup_{0 < t < T} \del H_t.
\endaligned
\ee
We impose a non-degeneracy condition on the averaged flux on the hypersurfaces.

\begin{definition}
\label{hyperb-def}
Let $M$ be a manifold endowed with a foliation \eqref{foli}-\eqref{foli2} and let
$\om=\om(\ubar)$ be a flux field. Then, the conservation law \eqref{LR.1} on $M$ satisfies the
{\rm global hyperbolicity condition} if there exist constants $0 < \barc < \cbar$ such that, 
for every non-empty hypersurface $e \subset H_t$,
the integral $\int_e i^*  \del_u \om(0)$ is positive and
 the function $\var_e: \RR \to \RR$,
$$
\var_e(\ubar) := \dashint_e i^* \om(\ubar)
               = {\int_e i^* \om(\ubar) \over \int_e i^*  \del_u \om(0)},
\qquad
\ubar \in \RR
$$
satisfies
\be
\label{hyperb}
\barc \leq \del_u \var_e(\ubar) \leq \cbar,  \qquad \ub \in \RR.
\ee
\end{definition}

The function $\var_e$ represents the {averaged flux} along $e$.
From now, we assume that the conditions above are satisfied and we refer to $H_0$ as an initial hypersurface
and we prescribe an initial data $u_0 : H_0 \to \RR$ on this hypersurface.
We impose a boundary data $u_B$ on the submanifold $B$. We sometimes refer to $H_t$ as {spacelike hypersurfaces.}

Under the global hyperbolicity condition \eqref{foli}-\eqref{hyperb}, the initial and boundary value problem takes
the following form. The boundary condition \eqref{LR.2i} decomposes into an initial data
\be
\label{cond1}
u_{H_0} = u_0
\ee
and a boundary condition
\be
\label{cond2}
u \big|_B \in \Ecal_{U,\Om}( u_B).
\ee
Correspondingly, the condition in Definition~\ref{key92} reads
$$
\aligned
&\int_M \Big(  d \psi \wedge \Om(u)
+ \psi \, (d  \Om) (u) - \psi \, \del_u U(u) (d \om) (u) \Big)
\\
& + \int_B \psi_{|\del M} \, \big(  i^*\Om(u_B) + \del_u U(u_B)  \big(\gamma - i^*\om(u_B) \big) \big)
  + \int_{H_T} i^* \Omega(u_{H_T}) - \int_{H_0} i^* \Omega(u_0) \geq 0.
\endaligned
$$
 
\begin{definition}
A flux field $\om$ is {\rm geometry-compatible} if it is closed for each value of the parameter,
\be
\label{LR.3}
( d \om ) (\ubar)=0, \qquad \ubar \in \RR.
\ee
\end{definition}

This condition  ensures that constants are trivial solutions, 
a property shared by many models of fluid dynamics (such as the shallow water model). 
When \eqref{LR.3} holds, it follows from Definition~\ref{key63} that
every entropy flux field $\Om$      satisfies 
$(d\Om) (\ubar) =0$ (for all $\ubar \in \RR$)
and the entropy inequalities \eqref{LR.1i} for a solution $u : M \to \RR$ take the simpler form
\be
\label{LR.1i-simple}
d( \Om(u)) \leq 0.
\ee


\subsection{The spacetime finite volume method}
\label{finit}
 
We now assume that $M= [0,T] \times N$ is foliated by slices with {\sl compact} topology $N$,
and the initial data $u_0$ is bounded.
We      assume that the global hyperbolicity condition holds and the flux field $\om$ is geometry-compatible. 
Let $\TT^h = \bigcup_{K \in\TT^h} K$ be a {triangulation} of $M$, that is,
a collection of {cells (or elements),} determined as the images of polyhedra of $\R^{n+1}$,
satisfying:
\begin{itemize}
\item The boundary $\dK$ of an element $K$ is a piecewise smooth, $n$-manifold,
$\dK = \bigcup_{e\subset \dK} e$ and  contains exactly two 
spacelike faces  $\ekp$ and $\ekm$ and ``vertical'' elements
$$
e^0 \in \del^0 K:= \del K \setminus \big\{\ekp, \ekm\big\}.
$$
\item The intersection $K \cap K'$ of two distinct elements $K, K' \in \TT^h$
is either a common face of $K, K'$ or else a submanifold with dimension at most $(n-1)$.
\item The triangulation is compatible with the foliation in the sense that
there exist times $t_0= 0 < t_1 < \ldots < t_N = T$ such that all spacelike faces
are submanifolds of $H_n := H_{t_n}$ for some $n=0, \ldots, N$,
and determine a triangulation of the slices.
We denote by $\Tcal^h_0$ the set of all $K$ which
admit one face belonging to the initial hypersurface $H_0$.

\end{itemize}

We define the measure $|e|$ of a hypersurface $e \subset M$ by
\be
\label{e-def}
|e|:= \int_e i^* \del_u \om(0).
\ee
This quantity is positive if $e$ is sufficiently ``close'' to one of the hypersurfaces along
which we have the hyperbolicity condition \eqref{hyperb}.
Provided $|e| >0$ which is the case if $e$ is included in one of the slices of the foliation,
 we associate to $e$ the function $\var_e: \RR \to \RR$. 
The following hyperbolicity condition holds along the triangulation since the 
spacelike elements are included in the spacelike slices: 
\be
\label{HYP}
\barc \leq \del_u  \var_{e_K^\pm}(\ubar) \leq \cbar,  \qquad \quad K \in \TT^h.
\ee

Next, we introduce the finite volume method by averaging \eqref{LR.1}
over each element $K \in \TT^h$. Applying Stokes theorem with a smooth solution $u$
to \eqref{LR.1}, we get
$$
0 = \int_K d(\om(u)) = \int_{\dK} i^*\om(u).
$$
Decomposing the boundary $\dK$ into its parts $\ekp, \ekm$, and $\dKz$, we obtain
\be
\label{LM.3}
\int_{\ekp} i^* \om(u) - \int_{\ekm} i^* \om(u) + \sumK \int_{e^0} i^*\om(u) = 0.
\ee
Given the averaged values $\ukm$ along $\ekm$ and $\ukezm$ along $\ekz \in \dKz,$ we need
an approximation $\ukp$ of the solution $u$ along $\ekp$.
The second term in \eqref{LM.3} can be approximated by
$$
\int_{\ekm} i^* \om(u) \approx \int_{\ekm} i^* \om(\ukm) = |\ekm|\var_{\ekm}(\ukm)
$$
and the last term by
$\int_{\ekz}i^*\om(u) \approx \qq (\ukm,\ukezm)$, 
where the \emph{total discrete flux} $\qq:\RR^{2}\to \RR$ (i.e., a scalar-valued function) must be prescribed.

Finally, the proposed version of the {finite volume method} for the conservation law \eqref{LR.1} takes the form
\be
\label{LM.4}
\int_{\ekp} i^*\om(\ukp) = \int_{\ekm} i^*\om(\ukm) - \sumK  \qq (\ukm,\ukezm)
\ee
or, equivalently,
\be
\label{LM.5}
|\ekp| \vep (\ukp) =  |\ekm| \vem(\ukm) -  \sumK \qq (\ukm,\ukezm).
\ee

We assume that the functions $\qq$ satisfy the following properties for all $\ubar,\vbar \in \RR:$
\begin{itemize}

\item \emph{Consistency:}
	\be
	\label{LM.3.1}
	\qq(\ubar,\ubar) = \int_{e^0} i^*\om(\ubar).
	\ee
	
\item \emph{Conservation:}
	\be
	\label{LM.3.2}
	\qq(\vbar,\ubar) = -q_{K_{e^0},e^0}(\ubar,\vbar).
	\ee
	
\item \emph{Monotonicity:}
	\be
	\label{LM.3.3}
	\del_{\ubar} \qq(\ubar,\vbar) \geq 0, \qquad \del_{\vbar}\qq(\ubar,\vbar)  \leq 0.
	\ee	
\end{itemize}

We need to specify the discretization of the initial data
and define constant initial values $u_{K,0} = \ukm$ (for $K \in \Tcal^h_0$)
associated with $H_0$, by setting
\be
\label{datai}
\int_{e_K^-} i^* \omega(u_K^-):= \int_{e_K^-} i^* \omega(u_0), \qquad  \, e_K^- \subset H_0.
\ee
We also define a piecewise constant function $u^h : M \to \RR$ by,  for every element $K \in \TT^h$, 
\be
\label{uh}
u^h(x) = \ukm, \qquad x \in K.
\ee

We  introduce $N_K := \# \del^0 K$,
 the total number of ``vertical'' neighbors of an element $K \in \Tcal^h$, 
supposed to be uniformly bounded. 
We fix a finite family of local charts covering the manifold $M$, and assume that the parameter $h$ coincides with the largest diameter of faces $e_K^\pm$ of elements $K \in \Tcal^h$,
where the diameter is computed with the Euclidian metric in chosen local coordinates.

We also impose the Courant-Friedrich-Levy condition (for all $K \in \TT^h$) 
\be
\label{CFL}
{N_K \over |\ekp|} \max_{e^0 \in \del^0 K} \sup_u \Big| \int_{e^0} \del_u \omega(u) \Big|  < \inf_u \del_u\vep,
\ee
in which the supremum and infimum in $u$ are taken over the range of the initial data. Finally, we assume that the family of triangulations satisfy 
\be
\label{AS2}
\lim_{h \to 0} {\tau_{\max}^2 + h^2 \over \tau_{\min}} = 
\lim_{h \to 0}
{\tau_{\max}^2 \over h} = 0 
\ee 
where $\tau_{\max} := \max_i (t_{i+1} - t_i)$ and $\tau_{\min} := \min_i (t_{i+1} - t_i)$. 
For instance, these conditions are satisfied if $\tau_{\max}$, $\tau_{\min}$, and $h$ 
vanish at the same order.

Our main objective in this presentation  is establishing the convergence of the proposed finite volume schemes towards
an entropy solution.
Our analysis of the finite volume method will rely on a decomposition of \eqref{LM.5} into (essentially) one-dimensional schemes, 
a technique that goes back to Tadmor \cite{Tadmor}, Coquel and LeFloch \cite{CL1}, and 
Cockburn, Coquel, and LeFloch \cite{CCL1}. 

By applying Stokes theorem to \eqref{LR.3} with some $\ubar \in \RR$, we obtain 
$$
\aligned
0 & = \int_K d(\om(\ubar))=\int_{\dK} i^*\om(\ubar)
\\
  & = \int_{\ekp} i^* \om(\ubar)- \int_{\ekm} i^* \om(\ubar)
       + \sumK  \qq(\ubar,\ubar).
\endaligned
$$
Choosing $\ubar= u_K^-$, we deduce  
\be
\label{ident}
|\ekp| \vep(\ukm) = |\ekm| \vem(\ukm)
- \sumK \qq(u_K^-, u_K^-),
\ee
which can be combined with \eqref{LM.5}: 
$$
\aligned
\vep( \ukp)
& = \vep( \ukm) -  \sumK {1 \over |\ekp|} \Big( \qq(\ukm,\ukezm) - \qq(\ukm,\ukm) \Big)
\\
& =  \sumK \left(
 {1 \over N_K} \vep( \ukm) - {1 \over |\ekp|} \Big( \qq(\ukm,\ukezm) - \qq(\ukm,\ukm) \Big) \right).
\endaligned
$$

We introduce the intermediate values $\tu$:  
\be
\label{LM.7}
\vep(\tu) :=  \vep( \ukm) -  {N_K \over |\ekp|} \Big( \qq(\ukm,\ukezm) - \qq(\ukm,\ukm) \Big),
\ee
and thus arrive at the {convex decomposition}
\be
\label{CD}
\vep(\ukp)=  {1 \over N_K} \sumK \vep(\tu).
\ee

Given any entropy pair $(U,\Om)$ and hypersurface $e\subset M$ satisfying $|e| >0$ we introduce
the averaged entropy flux along $e$: 
$\veO(u):= \dashint_e i^* \Om(u)$.  

\begin{lemma}
\label{convx}
For every convex entropy flux $\Omega$ one has
\be
\label{key94}
\vep^\Omega(\ukp) \leq  {1 \over N_K}  \sumK \vep^\Omega(\tu).
\ee
\end{lemma}

In fact, the function $ \vep^{\Om}\circ ( \vep^{\om})^{-1}$ is convex.

\noindent{\bf Proof.}  It suffices to show the inequality for the entropy flux, and then average
this inequality over $e$. We need to check 
\be
\label{318}
\Omega(\ukp) \leq  {1 \over N_K}  \sumK \Omega(\tu), 
\ee
namely 
$$
\aligned
&  {1 \over N_K}  \sumK \big( \Omega(\tu) - \Omega(\ukp) \big)
\\
& =  {1 \over N_K}
\sumK \big( \omega(\ukp) - \omega(\tu) \big) \del_u U(\ukp)
  +
 {1 \over N_K}  \sumK D_{K, e^0},
\endaligned
$$
with
$$
D_{K,e^0} : = \int_0^1 \del_{uu}U(\ukp) \Big(
\omega( \tu + a (\ukp - \tu)) - \omega(\tu) \Big) \, (\ukp - \tu) \, da.
$$
In the right-hand side, the former term vanishes identically (see \eqref{LM.7}) and 
the latter term is non-negative, since $U(u)$ is convex and $\del_u \omega$ is positive.
$\Box$


\subsection{Discrete entropy estimates}
\label{discret}
 
From the decomposition \eqref{CD}, we derive the discrete entropy inequalities of interest.

\begin{lemma}[Entropy inequalities for the faces]
\label{lem-1}
For all convex entropy pair $(U,\Om)$ and all $K \in\Tcal^h$ and $\ekz\in \del^{0}K$, 
there exists numerical entropy flux functions
$\QQ : \RR^{2}\to\RR$ satisfying (for all $u,v \in \RR$):
\begin{itemize}
\item $\QQ$ is consistent with the entropy flux $\Om$:
\be
\label{CONST}
\QQ(u,u ) = \int_{e^0} i^* \Om(u).
\ee
\item Conservation property:
\be \label{CONSV}
\QQ(u,v ) = -{Q}_{K_{\ekz},\ekz}(v,u).
\ee
\item Discrete entropy inequality:
       \be \label{DEI}
       \aligned
       \vepO&(\tu) - \vepO(\ukm)
       + {N_K \over |\ekp|}  \Big( \QQ(\ukm,\ukezm) -\QQ(\ukm,\ukm) \Big) \leq 0.
	\endaligned
	\ee
\end{itemize}
\end{lemma}

\noindent{\bf Proof.}  
{\it Step 1.} For  $u,v\in \RR$ and $  e^0\in \dKz$, let us set 
$$
\Hkez(u,v) :=\vep( u) - {N_K \over |\ekp|}  \Big( \qq(u,v) - \qq(u,u) \Big) 
$$
and note that
$\Hkez(u,u) = \vep( u)$. 
We now check that $\Hkez$ satisfies 
\be
\frac{\del}{\del u} \Hkez(u,v)\ge 0,
\qquad
\frac{\del}{\del v}\Hkez(u,v) \ge 0.
\ee
The second property is
immediate by the monotonicity \eqref{LM.3.3}. 
For the first one, we recall the
CFL condition \eqref{CFL} and the monotonicity \eqref{LM.3.3}. From the definition of $\Hkez(u,v)$,
we have 
$$
\Hkez(u,u_{K_{\ekz}}) = \big( 1- \sum_{\ekz\in \dKz}  \alpha_{K,\ekz}\big) \vep (u) + \sum_{\ekz\in \dKz}  \alpha_{K,\ekz} \vep (u_{K_{\ekz}}),
$$
and 
$$
\alpha_{K,\ekz}:= {1 \over |\ekp|} \frac{\qq(u,u_{K_{\ekz}})- \qq(u,u)}{\vep(u)- \vep(u_{K_{\ekz}})}.
$$
This gives a convex combination of $\vep(u)$ and $ \vep(u_{K_{\ekz}})$. By \eqref{LM.3.3}  we have $\sum_{\ekz \in \dKz} \alpha_{K,\ekz} \geq 0$ and, with \eqref{CFL},
$$
\sum_{\ekz \in \dKz} \alpha_{K,\ekz}
\leq
\sum_{\ekz \in \dKz} {1 \over |\ekp|} \Big| \frac {\qq(u,u_{K_{\ekz}})- \qq(u,u)}{\vep(u) - \vep(u_{K_{\ekz}})}\Big|  \leq 1.
$$

\

{\it Step 2.} We will establish the entropy inequalities for Kruzkov's entropies $\Omegabf$.
Introduce the discrete version of Kruzkov's entropy flux
$$
\Qbf(u,v,c):= \qq(u\vee c,v \vee c)- \qq(u\wedge c,v\wedge c),
$$
where $a\vee b= \max (a,b)$ and $a\wedge b= \min (a,b)$. 
Note that $\QQ(u,v)$ satisfies the first two properties of the lemma
with the entropy flux replaced by the Kruzkov's family $\Omega=\Omegabf$ in \eqref{KRZ}.

First, we observe: 
\be \label{LM.H1}
    \begin{aligned}
    &\Hkez(u \vee c,v \vee c) - \Hkez(u \wedge c,v\wedge c)
    \\[5pt]
    &=\vep(u\vee c) - {N_K \over |\ekp|} \big( \qq(u\vee c,v\vee c)
    - \qq(u\vee c,u\vee c)\big) \\[5pt]
    & \quad - \Big( \vep(u\wedge c) -  {N_K \over |\ekp|} \big(\qq(u\wedge c,v\wedge c)
    - \qq (u\wedge c,u\wedge c)\big) \Big) \\[5pt]
      &= \var^{\Omega}_{\ekp}( u ,c) - {N_K \over |\ekp|} \Big( \Qbf(u,v,c) -    \Qbf(u,u,c)\Big),
   \end{aligned}
\ee
where 
$\vep(u\vee c) - \vep(u\wedge c)= \dashint_{\ekp}i^*\Omegabf(u,c) =
\var^{\Omega}_{\ekp}( u ,c)$. 

Second, we prove that for $u= \ukm$, $v= \ukezm$ and for any $c\in \RR$
\be
\label{LM.H2}
\Hkez(\ukm \vee c,\ukezm \vee c) - \Hkez(\ukm \wedge c,\ukezm\wedge c)
 \geq \var_{\ekp}^{\Omega}( \tu,c).
 \ee
Indeed, we have 
$$
  \begin{array}{l}
     H_{K,e^0}(u,v) \vee H_{K,e^0}(\lambda,\lambda) \leq
     H_{K,e^0}(u\vee\lambda,v\vee\lambda), \\[5pt]
     H_{K,e^0}(u,v) \wedge H_{K,e^0}(\lambda,\lambda) \geq
     H_{K,e^0}(u\wedge\lambda,v\wedge\lambda),
   \end{array}
$$
where  $\Hkez$ is monotone in both variables. Since $\vep$ is monotone, we have
$$
\aligned 
& \Hkez(\ukm \vee  c,\ukezm \vee c) - \Hkez(\ukm \wedge c,\ukezm \wedge c)
\\
    &  \geq  \Big| \Hkez(\ukm,\ukezm)  -  \Hkez(c,c) \Big|
    =  \Big|  \vep(\tu)   - \vep(c) \Big|
    \\
    &  =
    \sgn\big(\vep(\tu)- \vep(c)\big) \big(\vep(\tu)- \vep(c)\big)
 \\
  & = \sgn\big( \tu - c \big)\big(\vep(\tu)- \vep(c)\big)
 =
    \varbO( \tu,c).
    \endaligned 
$$
Combining this with \eqref{LM.H1} (with  $u=\ukm$, $v=\ukezm $), we obtain the inequality 
$$
\aligned
\var_{\ekp}^{\Omega}&(\tu,c) - \var_{\ekp}^{\Omega}(\ukm,c)
     + {N_K \over |\ekp|} \Big( \Qbf(u,v,c) -    \Qbf(u,u,c)\Big) \leq 0, 
\endaligned
$$
which implies a similar inequality for all convex entropy flux fields. 
$\Box$


\

We now combine Lemma ~\ref{convx} with Lemma~\ref{lem-1}.

\begin{lemma}[Entropy inequalities for the elements]
For each $K \in \Tcal^h$, one has
\be
\label{DEI4}
|e_K^+| \, \big( \vepO (\ukp) - \vepO(\ukm) \big)
+ \sumK \big( Q(\ukm,\ukezm) -Q(\ukm,\ukm) \big) \leq 0.
\ee
\end{lemma}

If $V$ is convex, then a \emph{modulus of convexity} for $V$ is a positive real
$\beta < \inf V''$ (where the infimum is taken over the range of the data and solutions).
In view of the proof of Lemma \ref{convx}, $\var_e^{\Om}\circ (\var_e^{\om})^{-1}$ is convex for
every spacelike hypersurface $e$ and every convex function $U$.  (Note that 
the discrete entropy flux terms do not appear in \eqref{FVM.20} below.)

\begin{lemma}[Entropy balance inequality between two hypersurfaces]
\label{POC.i}
For $K \in\Tcal^h$, denote by $\beta_{e_K^{+}}$ a modulus of convexity for $\vepO \circ \big(\vepo\big)^{-1}$ and set $\beta = \min_{K \in \TT^h} \beta_{\ekp}$. Then, for $i \leq j$ one has
\label{PC-3}
\be
\label{FVM.20}
\aligned
& \sum_{K\in\Tcal_{t_j}^h} |\ekp| \vepO(\ukp)
  + \sum_{\substack{K\in\Tcal^h_{[t_i, t_j)}\\ \ekz\in\dKz}} \frac{\beta}{2 N_K} |\ekp| \big| \tu - \ukp \big|^2
 \leq \sum_{K\in\Tcal_{t_i}^h} |\ekm| \vemO(\ukm),
\endaligned
\ee
where $\Tcal^h_{t_i}$ is the subset of all $K$ satisfying $e_K^- \in H_{t_i}$, and  one
sets $\Tcal^h_{[t_i, t_j)} := \bigcup_{i \leq k < j} \Tcal^h_{t_k}$.
\end{lemma}

\noindent{\bf Proof.}  Multiplying \eqref{DEI} by $|\ekp|/N_K$ and summing in $K\in\Tcal^h$, $\ekz\in\dKz$ yield 
$$
\aligned
 \sum_{\substack{K\in\Tcal^h\\ \ekz\in\dKz}}
{|\ekp| \over N_K} \vepO (\tu ) 
& -  \sum_{K\in\Tcal^h}|\ekp|\vepO(\ukm)
\\
& 
 +  \sum_{\substack{K\in\Tcal^h\\ \ekz\in\dKz}} \big(\QQ(\ukm,\ukezm ) -\QQ(\ukm,\ukm) \big)
\leq 0.
\endaligned
$$
The conservation property \eqref{CONSV} gives
\be
\label{FVM.22}
  \sum_{\substack{K\in\Tcal^h\\ \ekz\in\dKz}}  \QQ(\ukm,\ukezm ) = 0
\ee
and so 
\be
\label{FVM.21i}
\aligned
&  \sum_{\substack{K\in\Tcal^h\\ \ekz\in\dKz}} {|\ekp| \over N_K} \vepO (\tu )  -  \sum_{K\in\Tcal^h}|\ekp|\vepO(\ukm)
 -  \sum_{\substack{K\in\Tcal^h\\ \ekz\in\dKz}} \QQ(\ukm,\ukm)
\leq 0.
\endaligned
\ee

If $V$ is convex and if $v =  \sum_{j}\alpha_{j} v_{j}$ is a convex combination of $v_{j}$,
then 
$$
V(v) + \frac{\beta}{2} \sum_{j}\alpha_{j} |v_{j}- v|^{2} \le \sum_{j} \alpha_{j} V(v_{j}),
$$
where $\beta =\inf V''$, the infimum being taken over all $v_j$. We apply this with  $ v= \vep(\ukp)$
and
$V = \vep ^{\Om}\circ (\vep^\omega)^{-1}$, which is convex.

In view of \eqref{CD} and
by multiplying the above inequality by $|\ekp|$ and summing in $K\in\Tcal^h$, we obtain
$$
\aligned
&  \sumkez |\ekp| \vepO(\ukp) +
  \sumkez {\beta \over 2}\, {|\ekp| \over N_K} \, |\tu - \ukp|^2
 \leq \sumkez {|\ekp|\over N_K} \vepO(\tu).
    \endaligned
$$
Combining the result with \eqref{FVM.21i}, we conclude that 
\be
\label{FVM.23}
\aligned
& \sum_{K\in\Tcal^h} |\ekp| \vepO(\ukp) - \sum_{K \in \Tcal^h} |\ekp| \vepO(\ukm)
 + \sumkez {\beta \over 2} {|\ekp| \over N_K} |\tu - \ukp|^2
\\
&  \leq  \sumkez  \QQ(\ukm,\ukm).
\endaligned
\ee
Finally, using 
$$
\aligned
 0 &=\int_{K} d(\Om(\ukm))= \int_{\del K} i^* \Om(\ukm)
\\
& = |\ekp| \vepO(\ukm) -  |\ekm| \vemO(\ukm)+ \sum_{\ekz\in\dKz} \QQ(\ukm,\ukm),
\endaligned
$$
we obtain the desired inequality, after further summation over all of $K$ within two arbitrary hypersurfaces. 
$\Box$

We apply Lemma~\ref{POC.i} and obtain an important uniform estimate.

\begin{lemma}[Global entropy dissipation estimate]
The entropy dissipation is globally bounded, as follows:
\be
\label{EDE}
\sum_{\substack{K\in\Tcal^h\\ \ekz\in\dKz}} \frac{|\ekp|}{N_K}\big| \tu - \ukp \big| ^2
\lesssim C \,  \int_{H_0} i^* \Omega(u_0)
\ee
for some constant $C>0$ depending upon the flux field and the sup-norm of the initial data. Here, 
$\Omega$ is the $n$-form entropy flux field associated with $U(u) = u^2/2$.
\end{lemma}

\noindent{\bf Proof.}  We apply \eqref{FVM.20} with the choice $U(u)= u^2$
$$
\aligned
0&\geq \sum_{K\in\Tcal^h} ( |\ekp| \vepO(\ukp) -|\ekm| \vemO(\ukm) )
 +  \sum_{\substack{K\in\Tcal^h\\ \ekz\in\dKz}} \
 \frac{\beta}{2} {|\ekp| \over N_K} \big| \tu - \ukp \big|^2.
\endaligned
$$
After summing up in the ``vertical'' direction and keeping the contribution of all 
$K \in \Tcal^h_0$ on $H_0$,
we deduce that 
$$
\aligned
  \sum_{\substack{K\in\Tcal^h\\ \ekz\in\dKz}}
  {|\ekp| \over N_K} \beta \big| \tu - \ukp \big|^2
  \leq
{2 \over \beta} \, \sum_{K\in\Tcal_0^h} |\ekm| \vemO(u_{K,0}).
\endaligned
$$
For some constant $C>0$, we have 
$\sum_{K\in\Tcal_0^h} |\ekm| \vemO(u_{K,0}) \leq C \, \int_{H_0} i^* \Omega(u_0)$.
These are essentially $L^2$ norm of the initial data, and this inequality is checked by fixing a reference volume form on $H_0$ 
and using the discretization \eqref{datai} of the initial data $u_0$. 
$\Box$


\subsection{Global form of the discrete entropy inequalities}
 
One additional notation now is needed in order to handle ``vertical face'' of the triangulation: we 
fix a reference field of non-degenerate $n$-forms $\omegab$ on $M$ (to measure the ``area'' 
of the faces $e^0 \in \del K^0$). This is used in the convergence proof only, but not in the formulation 
of the finite volume schemes. For every $K \in \Tcal^h$ we define 
\be
\label{e-def2}
|e^0|_{\omegab} := \int_{e^0} i^*\omegab \quad \text{ for faces } e^0 \in \del^0 K
\ee
and the non-degeneracy condition is equivalent to $|e^0|_{\omegab} >0$.  Given a smooth function $\psi$ defined on $M$ and given a face $e^0 \in \del^0 K$ of some element,
 we introduce 
$$
\psi_{\ekz} := {\int_{e^0} \psi \, i^* \omegab \over \int_{e^0} i^* \omegab}, 
\qquad
\psi_{\dKz} := {1 \over N_K} \sum_{\ekz\in\dKz}  \psi_{e^0}. 
$$ 

\begin{lemma}[Global form of the discrete entropy inequalities]
\label{PC-5}
Let $\Om$ be a convex entropy flux field and let $\psi \geq 0$ be a smooth function
supported away from the hypersurface $t=T$.
Then, the finite volume scheme satisfies the entropy inequality
\be
\label{FVM.31}
\aligned
& - \sum_{K \in \TT^h} \int_K d(\psi \Om ) (\ukm)
    - \sum_{K \in\Tcal^h_0} \int_{\ekm} \psi \, i^* \Om(u_{K,0})
\leq A^h(\psi) + B^h(\psi) + C^h(\psi), 
\endaligned 
\ee
with 
$$
\aligned 
& A^h(\psi) := \sumkez {|\ekp| \over N_K} \big( \psi_{\dKz} -  \psi_{\ekz} \big) \, \big( \vepO(\tu) - \vepO(\ukp) \big), 
\\
& B^h(\psi) := \sumkez \int_{\ekz} \big( \psi_{\ekz} - \psi \big) \, i^* \Om(\ukm), 
\\
& C^h(\psi) := -\sum_{K \in\Tcal^h} \int_{\ekp} \big( \psi_{\dKz} - \psi \big) \, \big( i^* \Om(\ukp) - i^*\Om(\ukm) \big).
\endaligned
$$
\end{lemma}

\noindent{\bf Proof.} 
From the discrete entropy inequalities \eqref{DEI}, we get 
\be
\label{FVM.36}
\begin{aligned}
& \sumkez {|e_K^+| \over N_K} \psi_{\ekz} \Big( \vepO (\tu ) - \vepO(\ukm)\Big)
\\
& \quad + \sumkez  \psi_{e^0}\ \big( \QQ(\ukm,\ukezm ) - \QQ(\ukm,\ukm) \big)
 \leq 0.
\end{aligned}
\ee
Thanks \eqref{CONSV}, we have
$\sumkez \psi_{\ekz} \QQ(\ukm,\ukezm ) = 0$
and, from \eqref{CONST}, 
$$
\aligned
&  \sumkez  \psi_{\ekz}\QQ(\ukm,\ukm ) =  \sumkez \psi_{\ekz} \int_{\ekz}i^*\Om(\ukm)
\\
& \quad =  \sumkez \int_{\ekz} \psi \, i^* \Om(\ukm)
    +  \sumkez \int_{\ekz} (\psi_{\ekz} - \psi)  i^* \Om(\ukm).
\endaligned
$$
Next, we observe 
$$
\aligned
& \sumkez {|\ekp| \over N_K} \psi_{e^0} \, \vepO(\tu)
\\
& = \sumkez {|\ekp| \over N_K} \psi_{\dKz} \vepO(\tu)
  +\sumkez {|\ekp| \over N_K} (\psi_{\ekz} - \psi_{\dKz} )  \vepO (\tu)
\\
&\geq \sum_{K \in\Tcal^h} |\ekp| \psi_{\dKz} \vepO(\ukp)
      +\sumkez  {|\ekp| \over N_K} (\psi_{\ekz} - \psi_{\dKz} ) \vepO ( \tu ), 
\endaligned
$$ 
where, we recalled \eqref{key94} and the convex combination \eqref{CD}.
From 
$$
\sumkez {|\ekp| \over N_K} \psi_{\ekz} \, \vepO(  \ukm)
=   \sum_{K\in\Tcal^h}|\ekp| \psi_{\dKz} \, \vepO(  \ukm), 
$$
the inequality \eqref{FVM.36} reads 
\be
\label{FVM.39}
\aligned
& \sum_{K \in\Tcal^h} |\ekp|\psi_{\dKz}  \Big(
 \vepO(\ukp) - \vepO(\ukm)\Big)
 -\sumkez \int_{\ekz} \psi \,  i^*\Om (\ukm)
\\
& \leq   -\sumkez {|\ekp| \over N_K} (\psi_{\ekz} - \psi_{\dKz}) \vepO ( \tu )
 + \sumkez \int_{\ekz} (\psi_{\ekz} - \psi)\, i^* \Om(\ukm).
\endaligned
\ee

The first term in \eqref{FVM.39} reads 
$$
\aligned
&\sum_{K \in\Tcal^h}|\ekp| \psi_{\dKz}  \Big(
\vepO(\ukp) -\vepO (\ukm)\Big)
 \\
 &= \sum_{K \in\Tcal^h} \int_{\ekp} \psi (i^* \Om(\ukp) - i^*\Om(\ukm))
  + \sum_{K \in\Tcal^h} \int_{\ekp} (\psi_{\dKz} - \psi) \,   (i^* \Om(\ukp) - i^*\Om(\ukm)).
\endaligned
$$
We sum up with respect to $K$ the identities 
$$
\aligned
 \int_{K} d(\psi \Om )(\ukm) 
& = \int_{\del K} \psi \, i^*\Om(\ukm)
\\
& =\int_{\ekp} \psi \, i^*\Om(\ukm) - \int_{\ekm} \psi \,i^*\Om(\ukm)
 + \sum_{\ekz\in\dKz} \int_{\ekz} \psi \, i^* \Om(\ukm)
\endaligned
$$
and we combine them with \eqref{FVM.39}. We arrive at the desired conclusion 
by observing that
$$
\aligned
& \sum_{K \in\Tcal^h} \Big(  \int_{\ekp} \psi \, i^* \Om(\ukp) -   \int_{\ekm} \psi \, i^* \Om(\ukm)\Big)
  = - \sum_{K \in\Tcal^h_0} \int_{\ekm} \psi \, i^* \Om(u_{K,0}).
\endaligned
$$
$\Box$


\subsection{Convergence and well-posedness results}
\label{wellp}
 
This is the final step of our analysis. 

\begin{theorem}[Convergence theory]
\label{convergencefvm}
Under the assum\-p\-tions in Section~\ref{finit}, 
the family of approximate solutions $u^h$ generated by the finite volume scheme
converges (as $h \to 0$) to an entropy solution to the initial value problem \eqref{LR.1}, \eqref{cond1}.
\end{theorem}

This theorem generalizes to spacetimes the technique originally introduced by
Cockburn, Coquel and LeFloch \cite{CCL00,CCL2} for the (flat) Euclidean setting and 
extended to Riemannian manifolds by \cite{ABL} and to Lorentzian manifolds by \cite{ALO}.

\begin{corollary}[Well-posedness theory on a spacetime]
\label{main}
Fix $\MM=[0,T] \times N$ a $(n+1)$-dimensional spacetime foliated by $n$-dimen\-sional
hypersurfaces $H_t$ ($t \in [0,T]$) with compact topology $N$ (cf.~\eqref{LR.1}).
Consider also  a geometry-compatible flux field $\om$ 
on $\MM$ satisfying the global hyperbolicity condition \eqref{hyperb}. 
Given any initial data $u_0$ on $H_0$,
the initial value problem \eqref{LR.1}, \eqref{cond1} admits
a unique entropy solution $u \in L^\infty(M)$ which has well-defined $L^1$ traces on spacelike hypersurface of $M$.
These solutions determines a (Lipschitz continuous) contracting semi-group: 
\be
\label{MTHM.1}
\int_{H'} i_{H'}^* \Omegabf\big( u_{H'}, v_{H'} \big)
\leq \int_{H} i^*_H \Omegabf\big( u_H,v_H\big)
\ee
 for any two hypersurfaces $H,H'$ such that $H'$ lies in the future of $H$, and 
the initial condition is assumed in the sense
\be
\label{MTHM.23}
\lim_{t \to 0 \atop t>0} \int_{H_t} \iht \Omegabf\big( u(t), v(t) \big) = \int_{H_0} \ihz \Omegabf(u_0,v_0).
\ee
\end{corollary}

The following conclusion was originally established by DiPerna \cite{DiPerna} for conservation laws posed on the Euclidian space. 

\begin{theorem}
\label{MTHM.2}
Fix $\om$ a geometry-compatible flux field on $\MM$ satisfying the global hyperbolicity condition \eqref{hyperb}.
Then, any entropy measure-valued solution $\nu$ 
to the initial value problem \eqref{LR.1}, \eqref{cond1} reduces to a Dirac mass and, more precisely, 
\be
\nu= \delta_u,
\ee
where $u \in L^\infty(\MM)$ is the entropy solution to the problem.
\end{theorem}

We now give a proof of Theorem~\ref{convergencefvm}.
By definition, a Young measure $\nu$ represents all weak-$*$ limits of composite functions $a(u^h)$
for all continuous functions $a$ (as $h \to 0$):  
\be
\label{COP.1}
a(u^h) \mathrel{\mathop{\rightharpoonup}\limits^{*}} \langle \nu, a \rangle
=
 \int_\RR a(\lambda) \, d\nu(\lambda). 
\ee

\begin{lemma}[Entropy inequalities for Young measures]
\label{COP-1}
Given any Young measure $\nu$ associated with the approximations $u^h$.
and for all convex entropy flux field $\Om$ and smooth functions $\psi \geq 0$ 
supported away from the hypersurface $t=T$, one has
\be
\label{COP.2}
\int_{M} \langle  \nu, d\psi \wedge \Om(\cdot)\rangle
- \int_{H_{0}} i^*\Om(u_0) \leq 0.
\ee
\end{lemma}
 
Thanks to \eqref{COP.2}, for all convex entropy pairs $(U,\Om)$ we have $d  \langle \nu, \Om(\cdot) \rangle  \leq 0$
  on $M$.
On the initial hypersurface $H_0$ the Young measure $\nu$ coincides with
the Dirac mass $\delta_{u_0}$.
By Theorem~\ref{MTHM.2} there exists a unique function $u \in L^\infty(M)$ 
such that the measure $\nu$ coincides with the Dirac mass $\delta_u$.
This  implies that $u^h$ converge strongly to $u$,
and this concludes our proof of convergence.

\noindent{\bf Proof.} We pass to the limit in \eqref{FVM.31},
by using the property \eqref{COP.1} of the Young measure.
Observe that  the left-hand side of \eqref{FVM.31} converges  to
the left-hand side of \eqref{COP.2}.
Indeed, since $\omega$ is geometry-compatible, the first term 
$$
\sum_{K \in \TT^h} \int_K d(\psi \Om ) (\ukm)
= 
\sum_{K \in \TT^h} \int_K d\psi \wedge \Om (\ukm)
= 
\int_M d\psi \wedge \Om(u^h) 
$$
converges to $\int_{M} \langle  \nu, d\psi \wedge \Om(\cdot)\rangle$. On the other hand,  
one has 
$$
\sum_{K \in\Tcal^h_0} \int_{\ekm} \psi \, i^* \Om(u_{K,0})
=
\int_{H_0} \psi \, i^*\Om(u_0^h) 
\to \int_{H_{0}} \psi \, i^*\Om(u_0),  
$$ 
in which $u_0^h$ is the initial discretization of the data $u_0$ converges strongly to $u_0$
since the maximal diameter $h$  tends to zero. 
  
   The terms on the right-hand side of \eqref{FVM.31} also vanish in the limit $h \to 0$.
We begin with the first term $A^h(\psi)$.
Taking the modulus, applying Cauchy-Schwarz inequality, and using \eqref{EDE}, we obtain 
$$
\aligned
|A^h(\psi)|  
&
\leq
\sumkez {|\ekp| \over N_K} |\psi_{\dKz} -\psi| |\tu - \ukm|
\\
&
\leq
\Big(\sumkez {|\ekp|\over N_K} |\psi_{\dKz} -\psi|^2\Big)^{1/2} \Big( \sumkez {|\ekp| \over N_K} |\tu - \ukm| ^2 \Big)^{1/2}
\\
& 
\leq \Big( \sumkez {|\ekp|\over N_K}  (C \, (\tau_{\max} + h))^2\Big)^{1/2} \Big( \int_{H_0}i^*\Om(u_0)\Big) ^{1/2}, 
\endaligned
$$
hence
$$
\aligned
|A^h(\psi) | 
& \leq C' \, (\tau_{\max} + h) \Big( \sumk |\ekp| \Big)^{1/2}
\leq C'' \,  {\tau_{\max} + h  \over (\tau_{\min})^{1/2}}. 
\endaligned
$$
Here, $\Om$ is associated with the quadratic entropy and 
we used that $|\psi_{\dKz} - \psi| \leq C \, (\tau_{\max} + h)$. 
Our conditions \eqref{AS2} imply that the upper bound for $A^h(\psi)$ tends to zero with $h$.

Next, we rely on the regularity of $\psi$ and $\Om$ and we estimate the second term in the right-hand side of \eqref{FVM.31}. 
By setting  
$C_{\ekz} := {\int_{\ekz} i^* \Omega(u_K^-) \over  \int_{\ekz} i^* \omegab}$, 
we obtain 
$$
\aligned
|B^h(\psi)| 
& = \Big| \sumkez \int_{\ekz} (\psi_{\ekz} - \psi)\, \Big( i^* \Om(\ukm) - C_{\ekz} i^* \omegab \Big) \Big| 
\\
& \leq \sumkez  \sup_K |\psi_{\ekz} - \psi| \, \int_{\ekz} \Big| i^* \Om(\ukm) - C_{\ekz} i^* \omegab \Big|
\leq C \, \sumkez  (\tau_{\max} + h)^2 \, |e^0|_{\omegab},  
\endaligned
$$
hence
$|B^h(\psi)| 
\leq C \, {(\tau_{\max} + h)^2  \over h}$. 
This implies the upper bound for $B^h(\psi)$ tends to zero with $h$. 

Finally, we treat the last term in the right-hand side of \eqref{FVM.31}
$$
\aligned
|C^h(\psi)| 
& \leq 
\sum_{K \in\Tcal^h} |\ekp| \sup_K |\psi_{\dKz} -\psi| \int_{\ekp}  |i^* \Om(\ukp) - i^*\Om(\ukm)|,  
\endaligned
$$
using the modulus defined earlier. 
In view of \eqref{318}, we obtain
$$
\aligned
|C^h(\psi)| 
& \leq
C \,  \sumkez {|\ekp| \over N_K} |\psi_{\dKz} -\psi| \,  \big| \tu - \ukm \big|, 
\endaligned
$$
and it is now clear that $C^h(\psi)$ satisfies the same estimate as the one we derived for $A^h(\psi)$.  
$\Box$


\subsection*{Acknowledgements}

The author was partially supported by the Agence Nationale de la Recherche (ANR) through the grant
ANR SIMI-1-003-01, and by the Centre National de la Recherche Scientifique (CNRS). These notes were written at the occasion of a short course given by the authors at the University of Malaga for the XIV Spanish-French School Jacques-Louis Lions. The author is particularly grateful to C. V\'azquez-Cend\'on and C. Par\'es for their invitation, warm welcome, and efficient organization during his stay in Malaga.

%

\end{document}